\documentclass[12pt]{preprint} 
\usepackage{times,amssymb,amsmath} 
\usepackage{mhenvs}      
\begin{document} 
\title{An $L^2$ theory for differential forms on path spaces I} 
\author{K.D. Elworthy and Xue-Mei Li} 
\date{}
 \institute{ 
Mathematics Institute, The University of Warwick \and
 Mathematical Sciences,
Loughborough University\\
\email{kde@maths.warwick.ac.uk}\\
\email{xue-mei.li@lboro.ac.uk}
} 
\maketitle 
 
\flushbottom

\newcommand{\C}{{\mathcal C}} 
\newcommand{\D}{{\rm I \! D}} 
\newcommand{\pathR}{{\mathcal{\rm I\!R}}} 
\newcommand{\Nabla}{{\bf \nabla}} 
\newcommand{\E}{{\mathbf E}} 
\newcommand{\F}{{\mathcal F}} 
\newcommand{\G}{{\mathcal G}} 
\newcommand{\HH}{{{\mathcal H}} }
\def\H{{\mathcal {H}}}
\newcommand{\K}{{\mathcal K}}

\def\P{{{\mathbf P}}} 
\def\p{{{\mathbf P}}} 
\newcommand{\R}{{{\mathbf R} }}
\def\Nabla{{{\,\nabla\!\!\!\!\!\!\nabla}}} 
\newcommand{\RR}{{\mathcal R}} 
\newcommand{\1}{{\bf 1}} 
\newcommand{\I}{{\mathcal I}} 
\newcommand{\W}{{\mathcal W}} 
\newcommand{\Sharp}{{\#}}
\def\L{{{\mathcal {L}}}}
\def\LL{{\mathbf {L}}}
\def\X{{\Bbb {X} }}
\def\Y{{\Bbb {Y} }}

\def\Ric{{\mathop {\rm Ric} }}
\def\cf{\textit{c.f. }}
\def\eg{\textit{e.g. }}
\def\as{\textit{a.s. }}
\def\ie{\textit{i.e. }}
\def\esssup{{\mathop{\rm ess\; sup}}}
\def\proj{{\mathop{\rm proj}}}
\def\dim{{\mathop{\rm dim}}}
\def\div{{\mathop{\rm div}}}
\def\ker{{\mathop{\rm ker}}}
\def\Hess{{\mathop{\rm Hess}}}
\def\Image{{\mathop{\rm Image}}}
\def\Dom{{\mathop{\rm Dom}}}
\def\dom{\mathop {\rm {\cal D}om}}
\def\id{{\mathop {\rm Id}}}
\def\Cyl{{\mathop {\rm Cyl}}}
\def\Span{{\mathop {\rm Span}}}
\def\s.t.{{\mathop {\rm s.t.}}}
\def\trace{{\mathop{\rm trace}}}
\def\ev{{\mathop {\rm ev}}}
\def\constant{{\mathop {\rm const.}}}
\def\var{{\mathop{\rm Var}}}
\def\supp{{\mathop{\rm Supp}}}
\def\im{{\mathop{\rm Im}}}
\def\Hom{{\mathop{\rm Hom}}}

\def\paral{{/\kern-0.55ex/} }
\def\parals_#1{/\kern-0.55ex/_{\!#1}} 
\def\n#1{|\kern-0.24em|\kern-0.24em|#1|\kern-0.24em|\kern-0.24em|}

\renewcommand{\thefootnote}{} 
 
\footnote{Research, in part,  supported by EPSRC GR/NOO 845. Elworthy benefited from an EU grant ERB-FMRX-CT96-0075, Li  from a Royal society Leverhulme Trust Senior  Research Fellowshop, NSF  research grant DMS 0072387,  and support from the Alexander von Humboldt Foundation.} 
 
\pagestyle{myheadings} 

\begin{abstract} 
An $L^2$ theory of differential forms is proposed for the Banach manifold of continuous paths on Riemannian manifolds $M$ furnished with its  Brownian motion measure.   Differentiation must be restricted to certain Hilbert space directions, the $H$-tangent vectors. To obtain a closed exterior differential operator the relevant spaces of differential forms, the $H$-forms, are perturbed by the curvature of $M$. A Hodge decomposition is given for $L^2$ $H$-one-forms, and the structure of $H$-two -forms is described. The dual operator $d^*$ is analysed in terms of a natural connection on the $H$-tangent spaces. Malliavin calculus is a basic tool.
\end{abstract} 
\bigskip

 \section{Introduction}
 
{\bf Background.}  We are concerned with  the construction of an $L^2$ Hodge theory on path spaces with respect to a suitable reference measure and a collection of `admissible' vector fields.
Consider the space of continuous paths on a compact Riemannian 
manifold, over a fixed time interval $[0,T]$. Path spaces are 
Banach manifolds with the usual concepts of differentiable 
functions and differential forms, for example see Eells \cite{Eells-setting}, Eliasson \cite{Eliasson67}, Lang \cite {Lang-book}   .
They also have a natural measure, their \emph{ Brownian motion}, or \emph{Wiener} measure.  
 
From the works of Bismut \cite{Bismut-Durham}, L\'eandre 
\cite{Jones-Leandre91}, Driver \cite{Driver92} and others following 
pioneering work by L. Gross \cite{GrossJFA1} in the classical Wiener 
 space case, it seems the natural Sobolev differential calculus for 
functions on  path spaces using such measures is of differentiation in directions given by 
Hilbert spaces of tangent vectors at each point: essentially the 
 tangent vectors of finite energy. These are the so called Bismut 
 tangent spaces. 
The integration by parts formula given by
 Driver \cite{Driver92}, and subsequent results suggest that 
these notions will lead to a satisfactory, and useful, Malliavin type 
 calculus in this context. However the construction of differential 
form theory using Bismut tangent spaces leads to difficulties even 
at the level of the definition of exterior derivative. This is 
because of the lack of integrability of Bismut tangent `bundle': 
the Lie bracket of suitable Bismut tangent space valued vector fields 
does determine a vector field, but in the presence of curvature it 
no longer takes values in the Bismut tangent spaces.  Several ways 
of getting round this problem have been formulated, and carried 
out, especially by L\'eandre \cite{Leandre-cohomology96} 
\cite{Leandre-homogeneous} 
 \cite{Leandre-stochastic-cohomology} who gave analytical de Rham groups 
 and showed  that they agree with the singular cohomology of the spaces. 
See also \cite{Leandre-homology}. But we are not aware of any which 
 have led to an $L^2$ theory with Hodge-Kodaira Laplacian on our path spaces in the presence of curvature. 
In flat Wiener space the problem does not arise and the $L^2$ 
theory was defined and shown to be cohomologically trivial by 
Shigekawa  \cite{Shigekawa-Hodge} \cite{Shigekawa2003}. See also Mitoma \cite{Mitoma} and Arai-Mitoma\cite{Arai-Mitoma}. 
For Abstract Wiener manifolds, a class of infinite dimensional manifolds with an integrable  Hilbert bundle of admissible directions, see Piech \cite{Piech82}. 
 For $M$ a compact Lie group with bi-invariant metric the corresponding 
results were proved by Fang and Franchi \cite{Fang-Franchi97}, but using the 
Bismut tangent spaces obtained from the flat left invariant connection 
 on $M$ so the problem again is avoided. They also considered loop groups, \cite{Fang-Franchi97}.
 For work done on `sub-manifolds' of Wiener space see 
 Airault-van Biesen \cite{Airault-VanBiesen-codimension},van Biesen \cite{VanBiesen-divergence} and especially
 Kusuoka \cite{Kusuoka-cohomology} \cite{Kusuoka-forms}, 
 Kazumi-Shigekawa \cite{Kazumi-Shigekawa}. These submanifolds were constructed to replicate loop spaces over Riemannian manifolds, with their natural ``Brownian bridge" measures. 
For a general  survey see L\'eandre \cite{Leandre-survey99}, and for a more introductory article concentrating on the approach taken here, see \cite {Elworthy-Li-ICM}.

Let $M$ be a compact $C^\infty$ Riemannian manifold. For a fixed 
positive number $T$, consider the space $C_{x_0}M$ of continuous 
paths $\sigma : [0,T]\to M$ starting at a given point $x_0$ of $M$, 
furnished with its natural structure as a $C^\infty $ Banach 
manifold and Brownian motion measure $\mu _{x_0}$. For smooth differential forms there are the de Rham cohomology groups 
 $H^q_{deRham}(C_{x_0}M)$. 
C. J. Atkin informs us that the techniques of \cite {Atkin-I, Atkin-II} can be extended to show that the de Rham groups would be 
equal to the singular cohomology groups, even though $C_{x_0}M$ 
does not admit smooth 
 partitions of unity,  and so trivial for $q\ge 0$ since based path 
spaces are contractible.  For related work, also see  Lempert-Zhang \cite{Lempert-Zhang}  on Dolbeault cohomology of a loop space. Since our primary interest is in the 
differential analysis associated with the Brownian motion measure 
 $\mu $ on $C_{x_0}M$, which could equally well be considered on  H\"older paths of any exponent smaller than a half, we could use 
 H\"older  rather than continuous paths and it is really only for notational convenience that we do not. In that case we would have smooth partitions of unity, see Bonic, Frampton \& Tromba \cite{Bonic-Framphton-Tromba}. However contractibility need not imply triviality of  the de Rham cohomology groups when some restriction is put on the spaces  of forms. For example if  $f:\R\to \R$ is given by $f(x)=x$ then $df$ determines a non-trivial 
class in the first bounded de Rham group of $\R$. If $f$ has value +1 for $x>1$ and $-1$ for $x<1$ then $df$ is non-trivial in $ L^2$ -cohomology.  In finite 
dimensions the $L^2$ cohomology of a cover $\tilde M$ of a compact 
manifold $M$ gives important topological invariants of $M$ even 
when $\tilde M$ is contractible, eg see Atiyah \cite{Atiyah76}; note also 
Bueler-Prokhorenkov \cite{Bueler-Prokhorenkov},   Ahmed-Stroock \cite{Ahmed-Stroock}, and Gong-Wang \cite{Gong-Wang}.

The Bismut tangent spaces $H^1_\sigma$ are defined by the parallel translation 
 $$\parals_t(\sigma):T_{x_0}M\to T_{\sigma(t)}M$$ of the Levi-Civita connection and consist of  those $v\in T_\sigma C_{x_0}M$ such that 
$v_t=\parals_t(\sigma)h_t$ for $h_\cdot \in L_0^{2,1}\left([0,T];T_{x_0}M\right)$. To have a satisfying $L^2$ theory of differential forms on $C_{x_0}M$ the obvious choice would be to 
consider `H-forms' i.e. for 1-forms these would be $\phi$ with 
$\phi_\sigma\in (H^1_\sigma)^*$, $\sigma\in C_{x_0}M$, and this 
agrees with the natural $H$-derivative $d_\HH f$ for $f:C_{x_0}M\to 
\R$. For $L^2$ q-forms the obvious  choice would be $\phi$ with 
$\phi_\sigma\in\wedge^q(H_\sigma^1)^*$, using here the Hilbert 
space completion for the exterior product.
An $L^2$-de Rham theory would come from the complex of spaces of $L^2$ sections
\begin{equation}\label{complex} 
\dots\stackrel{\bar d}{\to}  L^2\Gamma\wedge^q(H^1_\sigma)^* 
\stackrel{\bar d}{\to} L^2\Gamma\wedge^{q+1}(H^1_\sigma)^* 
\stackrel{\bar d}{\to}\dots 
\end{equation} 
where $\bar d$ would be a closed operator obtained by closure from the usual exterior derivative:  for  $V^j$, $j=1$ to $q+1$  $C^1$ vector fields, and $\phi$ 
 a differentiable one-form:  
\begin{equation}\label{Palais} 
\begin{array}{ll} 
&d\phi\left(V^1\wedge\dots \wedge V^{q+1}\right)\\ 
&=\frac{1}{q+1}\sum_{i=1}^{q+1} (-1)^{i+1} L_{V^i}\left[\phi\left(V^1\wedge\dots 
\wedge \widehat{V^i} \wedge\dots \wedge V^{q+1}\right)\right]\\ 
&+\frac{1}{q+1}{\sum_{1\le i<j\le q+1} (-1)^{i+j}  } 
 \phi\left( [V^i,V^j]\wedge  V^1\wedge\dots 
  \widehat{V^i}\wedge \dots \widehat{V^j}\dots 
 \wedge V^{q+1}\right) 
\end{array} 
\end{equation} 
where $[V^i,V^j]$ is the Lie bracket and $ \widehat{V^j}$ 
means omission of the vector field $V^j$. 

 From this would come the  de Rham-Hodge-Kodaira Laplacians $\bar d\bar d^*+\bar d^* \bar d$ and an 
 associated Hodge decomposition. However the brackets $[V^i, V^j]$ of sections of $H^1_\cdot$ are not in general  sections of $H_\cdot^1$, 
e.g. see Cruzeiro-Malliavin \cite{Cruzeiro-Malliavin96}, 
Driver \cite{Driver-Lie-bracket}, see also \cite{Elworthy-Li-vector-fields}, and formula (\ref{Palais}), below, for $d$  does not make sense for $\phi_\sigma$ defined only on 
 $\wedge^q H_\sigma^1$, each $\sigma$, as mentioned earlier. 

\medskip 
 
 Our proposal is to replace the Hilbert spaces $\wedge^q H_\sigma^1$ 
in (\ref{complex}) by a family of different Hilbert spaces 
 ${\mathcal H}_\sigma^q, q=2,3,\dots$, continuously included in $\wedge^q T_\sigma C_{x_0}M$, 
 though keeping the exterior derivative a closure of the classical exterior derivative on smooth cylindrical forms.
 
In Elworthy-Li \cite{Elworthy-Li-Hodge-1}, for $q=1,2$, we identified a class of   Hilbert subspaces $\HH _\sigma ^q$, of the completed exterior 
 powers $\wedge ^qT_\sigma C_{x_0}$ of the tangent space $T_\sigma C_{x_0}$ 
 to $C_{x_0}M$ at  a path $\sigma $ which could be the basic building blocks 
 of an $L^2$ de Rham and  Hodge theory for $C_{x_0}M$. We described $\HH _\sigma ^2$ without proof, proved 
closability of exterior differentiation on corresponding $L^2$ 1-forms, 
defined a self-adjoint Hodge Kodaira Laplacian on such $L^2$ 1-forms and 
established the Hodge decomposition. 

The article \cite{Elworthy-Li-vector-fields} both discusses some of the constructions here for more general diffusion measures and connections, and relates them to the Bismut type formulae for differential forms on $M$, \cite{Elworthy-Li-forms-CR}, see also \cite Driver-Thalmaier {Driver-Thalmaier}. In particular it shows that a very natural class of two-vector fields on
$C_{x_0}M$ are of the type we consider here (\ie are sections of $\HH^2$). 
\medskip

{\bf Main Results.}
 Here we give a detailed analysis of $\HH _\sigma ^2$ and define $\HH^q_\sigma$ for $q>0$.  For $q=1$, as a space $\HH^1_\sigma=H_\sigma^1$. For flat manifolds,  $\HH^q_\sigma=\wedge ^q H^1_\sigma$ for all $q$ and the standard Hodge decomposition theorem follows.  However in general, the spaces $\HH^q_\sigma$ we construct are different from $\wedge ^q H$, the exterior products of the Bismut tangent bundle.
Sections of  $\HH^q$ are called $H$-q-vector fields and sections of  $(\HH^q)^*$ $H$-differential forms of degree $q$. In fact   $\HH^2_\sigma$ is a deformation of $\wedge^2 H^1_\sigma$ inside $\LL_{skew}(\HH_\sigma^1 ,\HH_\sigma^1)$ by the curvature of $M$. As a Hilbert space $\HH^2_\sigma$ is defined to be isometric to $\wedge^2 H^1_\sigma$ by a map involving the curvature of the so called damped Markovian connection on the Bismut tangent ``bundle".  
 Algebraic operations such as interior products acting on H-two vectors, and the exterior products of $H$-one forms, as well as the derivation property for the exterior derivative are shown to make sense.  A Hodge decomposition is given for H-one-forms.   In a sequel, Part II,  we  establish the analogous decomposition for $L^2$ $2$ -forms, and  we show that  the spaces  $\HH ^q_\sigma$  defined by suitable 
 It\^o maps $\I$  depend only on the Riemannian structure of the base manifold $M$.

 \medskip
 
 \noindent {\bf Organisation. }
 The article is organised as follows:\\
\begin{tabular}{ll} 
\S\ref{se-exterior}. \hskip 4pt & Review of basic results concerning exterior powers of 
relevant spaces \\ 
& of tangent vectors to $C_{x_0}M$. \\ 
 \S\ref{se-define-Hq}. & Special It\^{o} maps and the definition of $\HH  ^q$. \\ 
\S\ref{se-characterization}. & Characterisation of $\HH^1$ and $\HH ^2$. \\ 
\S \ref{se-1-form}. & H-one-forms: exterior differentiation and Hodge decomposition. \\ 
\S \ref{se-tensor}.  & Tensor products as operators: algebraic operations on H-one forms . \\ 
\S\ref{se-derivative-p}. & The derivation property of $\bar{d^1}$ . \\ 
\S \ref{se-divergences}. & Infinitesimal rotations as divergences. \\ 
 \S\ref{se-torsion}. & Differential geometry of the space $\HH^2$ of two-vectors. \\ 
 Appendix A. & Conventions.\\
 Appendix B. &  Brackets of vector fields, torsion, and $d\phi(v^1\wedge v^2)$.\\

\end{tabular} 
\medskip

In \S\ref{se-exterior} we discuss the various completed tensor 
 products of tangent, and other spaces which we will use. Properties of these relating to tensor products of abstract Wiener spaces are used in order to define our spaces $\HH^q$ in \S\ref{se-define-Hq}.   The aim is to show that these constructions are well behaved and have interesting geometry.
 \medskip
 
 \noindent  One of the main results, see  \S\ref{se-characterization},  is a characterisation of $\HH^2$ as a perturbation of $\wedge^2\HH$ by a curvature of the Levi-Civita connection on $M$.
 Write $\HH_\sigma=\HH^1_\sigma$,  then
\begin{eqnarray}
\HH^2_\sigma&=&(I+Q_\sigma)\wedge^2 \HH_\sigma
\end{eqnarray}
 for some operator $Q_\sigma$ on $\wedge^2_\epsilon T_\sigma C_{x_0}$. Equivalently 
 $$u\in \HH^2 \hbox {   if and only if } u-\pathR(u)\in \wedge^2\HH$$
 where $\pathR$ is identified in \S 9 as the curvature of the damped Markovian connection on the $H$-tangent spaces.
 \medskip
 
  In \S\ref{se-1-form} we rapidly recall the results concerning closability of our exterior derivative on H-one-forms and  the Hodge decomposition for H-one-forms.

  \medskip
 
 The remainder, the main part,  of the article is an analysis of the space $\HH^2$, its associated H-two-forms, and the adjoint of the exterior derivative, an operator from H-two-forms to H-one- forms, together with the corresponding divergence operator from two-vector-fields to vector fields.
   In \S\ref{se-tensor} it is shown that the exterior product of two H-one- forms is naturally an H-two form, and the interior product of an H-two form  with a $H$-one form is a $H$-one form. 
   The operator $Q$ has image in $\L_{skew}(\HH;\HH)$,
 which implies an element of $\HH^2_\sigma$  can be considered  
to be an element of $\L_{skew}(\HH_\sigma; \HH_\sigma)$, \cf Corollary \ref{H2-inclusion}, although in general it is not compact and so not in $\wedge^2\HH_\sigma$.  In \S\ref{se-derivative-p} a corresponding derivation formula for the exterior derivative of H-one-forms,  Theorem \ref{op-theorem-2}, is shown to hold.
  
  In \S \ref{se-divergences} it is shown that the elements of the image of  suitable smooth sections of $ \wedge^2\HH$ by $Q$
  ``have a divergence"  in the sense of satisfying an integration by parts formula and 
a formula is given in \ref{pr-divQ} for $\div Q(V^1\wedge V^2)$.
   Vector fields which are not H-vector fields also make their appearance, especially as Lie brackets. The latter involve infinitesimal rotations which ``have a divergence", and in their case the divergence is zero. It is natural to ask if they themselves are divergences, in this extended sense, of some two-vector field. In \S\ref{se-divergences} this is shown to be true in a wide class of adapted situations on flat Wiener space, Proposition \ref{pr-rot}. This has  independent interest, but it is extended, in Theorem \ref{th-divQ-2}, to show that the torsion of the damped Markovian connection when applied to suitable non-anticipating H-vector fields is the divergence of the perturbing factor in the definition of $\HH^2$:
\begin{equation} 
 \div Q(u^1\wedge u^2)={1\over 2}{\Bbb T}(u^1,u^2),
\end{equation} 
 Here ${\Bbb T}$ is the torsion of the damped Markovian connection $\Nabla$.
 This helps explain the ``cancellation" of the bracket occurring with our exterior derivative, and fits in with the result of Cruzeiro-Fang, \cite{Cruzeiro-Fang}, concerning the vanishing of the divergence of such torsions.
 The damped Markovian connection, introduced by Cruzeiro-Fang, \cite {Cruzeiro-Fang}, plays an important role here, as it did in \cite{Elworthy-Li-Ito-map}. As in \cite{Elworthy-Li-Ito-map} we introduce it by giving a $C_{id}([0,T];O(n))$-bundle structure to $\HH$. This is done in \S\ref{se-torsion}. Here we also relate the divergence of our H-two-vector fields to the adjoint of the damped Markovian covariant derivative in a non-anticipating situation, Corollary \ref{cor-nabla-star}:
  For suitable non-anticipating $U$, $V$, 
\begin{equation} 
{\Nabla}^*(U\wedge V)=\div(I+Q)(U\wedge V) 
\end{equation} 

 We also describe the curvature of the damped Markovian connection in \S\ref{se-torsion}D, to establish our claim that $\HH^2$ is a perturbation of $\wedge^2\HH^1$ using this curvature operator,  Theorem \ref{Theorem-3.2}(iii). In \S\ref{se-torsion}D we essentially show that $\D^{2,1}$ H-two-forms are in the domain of the adjoint of $\bar{d^1}^*$, extending the result for one-forms proved in \cite{Elworthy-Li-Ito-map}.

 \bigskip
 
\noindent 
{\it List of symbols.}\begin{itemize}\parskip=-0.15mm
\item[]  $C_{x_0}M$ or $C_{x_0}$ ---
 space of continuous paths over $M$ starting from $x_0$.\item[]   
$T_\sigma C_{x_0}M$, or $T_\sigma C_{x_0}$---  tangent space at $\sigma$ to $C_{x_0}M$.\item[]   
$\HH ^1_\sigma$ or $\HH_\sigma$--- Bismut tangent space, 
a Hilbert space included in $T_\sigma C_{x_0}$.\item[]   
$\HH ^1$ or $\HH$ ---  corresponding Bismut tangent ``bundle", $\cup \HH^1_\sigma$ \item[]   
$\HH ^2$ ---  vector ``bundle'' with fibres
   $\HH ^2_\sigma\subset \wedge^2 T_\sigma C_{x_0}M$.\item[]   
$\Gamma B$ ---  sections of a vector bundle $B$.\item[]   
$L^2\Gamma B$ ---  $L^2$ sections of a vector bundle $B$.\item[]   
$C_0\R^m$ --- Wiener space with Wiener measure $\P$, the canonical probability space.\item[]   
$L_0^{2,1}(G)$ --- For $G$ a Hilbert space, this is 
$\{ h: [0,T]\to G \hbox{ such  that } \int_0^T |\dot h_s|^2\, ds<\infty \}$. When $G=\R^m$, this is the Cameron-Martin space, denoted by  $H$.\item[]   
$(\xi_t, t\ge 0)$ --- a Brownian stochastic flow of diffeomorphisms of $M$ 
\item[]   
$T\xi_t$ --- space derivative of $\xi_t$.\item[]   
$\mu $ ---  Brownian motion measure, also called Wiener measure, on $C_{x_0}M$
\item[]   
$\I$ ---the It\^o map induced by $(\xi_t(x_0), t\ge 0)$, 
 $\I(\omega):=\xi_\cdot(x_0,\omega)$.\item[]   
 $T\I$ ---H-derivative of the It\^o map.\item[]   
$\F^{x_0}$ --- the algebra generated by   $(\xi_t(x_0), t\ge 0)$ on $M$.\item[]   
$\overline{f}(\sigma)$ ---  conditional expectation of $f$ given $\I=\sigma$, $\sigma\in C_{x_0}$, 
 e.g. $\overline{T\I}_\sigma$.\item[]   
$W^{(q)}_t$ --- Weitzenbock flow of $q$-vectors, equation (\ref{damped}).  \item[]   
${W^{(q)}}_t^s$ --- Weitzenbock flow starting from time $s$.\item[]   
$W_t$ ---   damped parallel translation, $W_t=W_t^{(1)}$. \item[]   
 
$\frac{\D^q}{dt}, \frac{\D}{dt}$---see \S\ref{se-characterization}.\item[] 
$L_2T_\sigma C_{x_0}$ --- the space of $L^2$ tangent vectors at $\sigma$,  Definition \ref{L2T}
 \item[]  
${\cal W}$ --- isometry  between $\HH^1$  and $L_2T C_{x_0}$, equation (\ref{w}).  
    \item[] 
$\L  (E_1;E_2)$ --- the space of continuous linear maps between  Banach spaces
\item[] $\L_2(H_1;H_2)$--- Hilbert -Schmidt maps between Hilbert spaces.

$\RR$, $\RR^q$, $Ric$ ---Respectively the curvature operator, the 
Weitzenbock curvature on $q$ forms, and the Ricci curvature on $M$.
\end{itemize}

\section{Exterior Powers; Notation} 
\label{se-exterior} For convenience the conventions we use for tensor products, exterior powers etc. are gathered together as an Appendix. Please note that they differ from those used in our previous articles, such as \cite{Elworthy-Li-Hodge-1}.
 
{\bf A.} All linear spaces are over $\R$. We shall deal with tensor products of  Hilbert spaces and of 
Banach spaces of continuous paths. For any linear space $E$ let 
$\otimes_0^q E$ denote the $q^{th}$-algebraic tensor product 
of $E$ with itself and ${\wedge_0^qE}$ the linear subspace of 
antisymmetric elements. For infinite dimensional Banach spaces $E$ we will need 
completions of these spaces, \eg see Ruston \cite{Ruston} or Cigler-Losert-Michor \cite{Cigler-Losert-Michor}: 
 
\begin{enumerate} 
\item[(i)]  When ${E=T_\sigma C_{x_0}}$ or $C_0 \R ^m$ let $%
\otimes ^qE$ and $\wedge ^qE$ denote the completions using the largest cross 
norm, i.e. the projective tensor products,
$$\|v\|_\pi=\inf\Big \{\sum_{i=1}^n \|a_i\|\|b_i\|, \hbox { where }
v=\sum_{i=1}^n a_i\otimes b_i, a_i, b_i\in E, n<\infty\Big\}.$$
\item[(ii)]  When $E$ is a Hilbert space $H$, let $\otimes ^qH$ and $\wedge 
^qH$ denote the standard Hilbert space completions, (so $\otimes ^2H$ can be 
identified with the space of Hilbert-Schmidt operators on $H$). 
 
\item[(iii)]  In general let $\otimes _\varepsilon ^qE$ and $\wedge 
_\varepsilon ^qE$ refer to the completions with respect to the smallest 
reasonable cross norm, i.e. the inductive cross norm,
$$\|w\|_\epsilon =\sup_{|u^*|\le 1, \|v^*\|\le 1, u^*,v^*\in E^*} \|u^*\otimes v^*(w)\|.$$
\end{enumerate} 

We shall use the natural inclusion maps as identifications and so consider $$\otimes_0^q E\subset \otimes ^qE\subset \otimes_\varepsilon^qE.$$ 
 
Thus a differential $q$-form $\phi $ on $C_{x_0}M$ which by definition gives a continuous antisymmetric multi-linear map $ \phi _\sigma 
: T_\sigma C_{x_0}\times \dots \times T_\sigma C_{x_0}\to R$, Lang\cite {Lang-book}, can equivalently be defined as a section of the bundle $ \L  \left( \wedge ^qTC_{x_0}; \R \right) $ with 
fibres the dual spaces ${(\wedge ^qT_\sigma C_{x_0})^{*},\sigma 
\in C_{x_0}M}$.

\bigskip 
 
\noindent 
{\bf B.} If $S:E_1\to E_2$ and 
$T:F_1\to F_2$ are two linear maps of linear spaces,  there is 
the induced linear map $S\otimes T: E_1\otimes_0 F_1\to E_2\otimes_0 F_2$. 
 The Banach space  constructions are functorial so that if ${S,T\in 
\L  \left( C_0 \R ^m;T_\sigma C_{x_0}\right) }$ then $S\otimes T$ 
determines a continuous linear map of the completed tensor spaces ${\otimes ^2C_0 \R ^m} 
$ to ${\otimes ^2T_\sigma C_{x_0}M}$ and if $S=T$ we have its restriction $%
{\wedge ^2S:\wedge ^2C_0 \R ^m\to \wedge ^2T_\sigma C_{x_0}M%
}$,   Ruston p63 \cite{Ruston} and Cigler-Losert-Michor \cite{Cigler-Losert-Michor}; with the corresponding 
 result for the inductive tensor  product, for the Hilbert space case, and for $q>2$.  There is also 
the estimate on operator norms 
\[ 
\|S^1\otimes \dots \otimes S^q\|\le \|S^1\|\cdot \dots \cdot \|S^q\| 
\] 
so that in particular 
\[ 
\|\wedge ^qS\|\le \|S\|^q 
\] 
in all of these cases, see Ruston \cite{Ruston} and Cigler-Losert-Michor \cite{Cigler-Losert-Michor}. \bigskip 
 
For example let $H\equiv L_0^{2,1} \R ^m$ be the (Cameron-Martin) 
 Hilbert space of functions ${h:[0,T]\to  \R ^m}$ of 
the form $h_t=\int_0^t  \dot h_s \, ds$ with $\dot h\in L^2([0,T]; \R ^m)$ and inner product 
${\langle h^1, h^2 \rangle 
=\int_0^T \langle \dot h_s^1, \dot h_s^2\rangle_{ \R ^m}\, ds}$. 
Thus the indefinite integral
$$\int_0^\cdot : L^2 \big([0, T];  \R ^m \big)\to H  $$ 
is an isometry with inverse which we will write as 
$${\frac{d}{d\cdot}}:  H \to L^2  \big([0,T]); \R ^m\big).  $$ 
From this we obtain the isometry 
$$\wedge^q \left(\int_0^\cdot\right): \wedge^q L^2 \big([0, T];  \R ^m 
\big)\to \wedge^q H   $$ 
with inverse 
$$\wedge^q\left({\frac{d}{d\cdot}}\right) : \wedge^q L_0^{2,1} \R %
^m\to \wedge^q L^2 \big([0,T]); \R ^m\big).  $$ 
 
\bigskip 
 
\noindent 
{\bf C. } We will regularly make use of the well known isometries 
\[ 
\otimes _\varepsilon ^qC_0 \R ^m\stackrel{\rho }{\longrightarrow }%
C_0\left( [0,T]^q;\otimes ^q \R ^m\right) 
\] 
where the right hand side consists of those continuous ${\alpha :[0,T]^q\to \otimes ^q \R ^m}$ for which ${\alpha 
(t_1,\dots ,t_q)=0}$ if $t_j=0$ for any $j$. For example see Cigler-Losert-Michor \cite{Cigler-Losert-Michor} p66. For $V\in \otimes_\varepsilon ^qC_0 \R ^m$, write 
\[ 
{V_{t_1,\dots ,t_q}:=\rho (V)(t_1,\dots ,t_q)}. 
\] 
Let ${{\rm ev}_t:C_0 \R %
^m\to  \R ^m}$ be the evaluation map at time $t$,  then 
$$V_{t_1,\dots ,t_q}=({\rm ev}_{t_1}\otimes 
\dots \otimes {\rm ev}_{t_q})V.$$
 
Also note that such $V$ lies in $\wedge _\varepsilon ^qC_0 \R ^m$ if and 
only if ${\rho (V):[0,T]^q\to \otimes ^q \R ^m}$ 
anti-commutes with permutations, i.e. 
\[ 
V_{t_{\pi (1)},\dots ,t_{\pi (q)}}=(-1)^{\pi }S_\pi V_{t_1,\dots ,t_q} 
\] 
for any permutation $\pi $ on $\{1,\dots ,q\}$ with $S_\pi $ the induced 
action on $\otimes ^q \R ^m$. If so 
\[ V_{t,\dots ,t}\in \wedge ^q \R ^m \] 
and 
\[ V_{t_1,t_2,\dots ,t_2}\in  \R ^m\otimes \wedge ^{q-1} \R ^m \] 
etc. From this we see that elements of $\wedge _\varepsilon ^qC_0 \R ^m$ 
and hence those of the smaller spaces $\wedge ^qC_0 \R ^m$ are determined by their values on the simplex $0\le t_1\le \dots \le t_q\le T$. 
 
Similarly, to any ${V\in \otimes ^q_\varepsilon T_\sigma C_{x_0}}$ we have $%
{V_{t_1.\dots ,t_q}\in T_{\sigma _{t_1}}M\otimes \dots \otimes 
T_{\sigma _{t_q}}M}$ corresponding to an isometric isomorphism of $\otimes _\varepsilon 
^qT_\sigma C_{x_0}$ with the space of continuous maps $V_{\cdot }$ on $%
[0,T]^q$ such that 
 
\begin{picture}(200, 100)(0,0) 
\put(215,76) {$\otimes^q TM$} 
\put(230,67){\vector(0,-1){40}} 
\put(80,12){\vector(1,0){115}} 
\put(40,10){$[0,T]^q$} 
\put(100,15){$\sigma_\cdot\times\dots \times\sigma_\cdot $} 
\put(67,27){\vector(3,1){140}} 
\put(200,10){$M\times \dots \times M$} 
\put(232, 40){$\pi$} 
\put(127, 57){$V$} 
\end{picture}\linebreak commutes and $V_{t_1,\dots, t_q}=0$ when $t_j=0$ for any $j$. 
 
\bigskip 
 
{\bf D.} By functorality the inclusion ${i:L_0^{2,1} \R %
^m\to C_0 \R ^m}$ gives rise to a continuous linear inclusion $%
{\otimes^q i:\otimes ^qH \to \otimes }%
_\varepsilon ^q{C_0 \R ^m}$. From \S B we see that ${V\in %
\hbox{Image}\otimes^q i}$ if and only if 
\begin{equation} 
    \label{completion-tensor-s} 
V_{t_1,\dots ,t_q}=\int_0^{t_1}\int_0^{t_2}\dots \int_0^{t_q 
                     }U_{s_1,\dots ,s_q}\,\, ds_1\dots ds_q, 
\end{equation} 
$(t_1,\dots ,t_q)\in [0,T]^q$, for some ${U_{\cdot }\in 
 L^2([0,T]^q;\otimes  \R ^m)}$. 
Here we use the isometry $\rho $ of ${\otimes ^qL^2([0,T]; 
 \R ^m)}$ with ${L^2([0,T]^q;\otimes ^q \R ^m)}$. In fact 
the $U_{\cdot }$ in the above formula is just ${\rho \left( 
\otimes ^q({\frac d{d\cdot }})V\right) }$ or equivalently 
 ${U_{t_1,\dots ,t_q}}$ is the weak derivative 
 ${\ {\frac{\partial ^q} 
  {\partial t_1\dots \partial t_q}}V_{t_1,\dots ,t_q}}$. 
 
\bigskip 
 
{\bf E.} Given a bounded linear operator $S:E\to F$ of Banach spaces there 
is also the functorial construction 
\[ 
\left( d\otimes ^q\right) (S)\,\,\,:\,\,\,\otimes _0^qE\longrightarrow 
\otimes _0^qF 
\] 
defined by linearity and 
\begin{eqnarray*} 
&&\left( \left( d\otimes ^q\right) (S)\right) (e^1\otimes \dots \otimes e^q) 
\\ 
&=&S(e^1)\otimes e^2\otimes \dots \otimes e^q+e^1\otimes S(e^2)\otimes \dots 
\otimes e^q+\dots +e^1\otimes e^2\otimes \dots \otimes S(e^q). 
\end{eqnarray*} 
This is just a sum of operators described in \S B and so extends over the 
relevant completion. The same notation will be kept for these extensions. 
\bigskip 
 
Note that if $V$ is in ${\otimes ^q H }$ then $\Big( \big( {d\otimes ^q}\big) ({\frac d{d\cdot }}) \Big) (V)$ is in $\otimes ^qL^2\big( [0,T]; \R %
^m\big)$ with kernel 
\begin{equation} 
\left( \left( {d\otimes ^q}\right) ({\frac d{d\cdot }})\right) 
(V)_{t_1,\dots ,t_q}=\sum_{j=1}^q{\frac \partial {\partial t_j}}V_{t_1,\dots 
,t_q}.  \label{section-1-1} 
\end{equation} 
The restriction $\left( d\Lambda ^q(S)\right) $ of $\left( d\otimes 
^q(S)\right) $ to $\wedge _0^qE$ has the form 
\begin{eqnarray*} 
&&\left( d\Lambda ^q(S)\right) (v^1\wedge v^2\wedge \dots \wedge v^q) \\ 
&=&S(v^1)\wedge v^2\wedge \dots \wedge v^q+\dots +v^1\wedge v^2\wedge \dots 
\wedge S(v^q) 
\end{eqnarray*} 
and for $q=2$ 
\[ 
\left( d\Lambda ^2(S)\right) (v^1\wedge v^2)={\frac 12}\left\{ Sv^1\otimes 
v^2+v^1\otimes Sv^2-Sv^2\otimes v^1-v^2\otimes Sv^1\right\} . 
\]

In general we shall use $|\;\;|$ to denote norms of finite dimensional 
spaces. $\|\;\|$ for infinite dimensional spaces, with 
$\n{ \;\;}$ for spaces such as 
 $L^2(\Omega ; \R ^n)$, or $L^2(C_{x_0}M; \R )$, where integration 
 over probability spaces are  involved. 
 
\section{Special It\^o maps and the definition of $\HH_\sigma^q$} 
\label{se-define-Hq} 
{\bf A.} Take a surjective $C^\infty$  vector bundle morphism, $X:\underline \R^m\to TM$, of the trivial $\R ^m$ bundle over $M$ onto $TM$, for some $m\ge n=\dim M$.  Suppose that  $X$ induces the given Riemannian metric on $M$ and let $Y$ be the $\R^m$ valued 1-form such that $Y_x=X(x)^*: T_xM\to \R^m$.
For $U$ a vector field and $v\in T_xM$, set \begin{equation}
\nabla_vU =X(x) d\left[ y\to Y_yU(y)\right](v),
\end{equation}
as  in Elworthy-LeJan-Li\cite{Elworthy-LeJan-Li-Tani}  \cite{Elworthy-LeJan-Li-book}, where it was called LW connection for $X$. 
Suppose that the connection $\nabla$ 
is the Levi-Civita connection.
Take $(B_t)$ to be the canonical Brownian motion on $ \R ^m$ 
 with probability space ${C}_0 \R ^m$ and Wiener measure 
 ${\Bbb P}$ and consider the stochastic differential equation on $M$
\begin{equation} 
dx_t=X(x_t)\circ dB_t,\hskip 40pt
0\le t\le T.
  \label{1} 
\end{equation} 
Then the  solutions are Brownian motions on $M$. Let $\mu _{x_0}$ be the Brownian motion measure on $C_{x_0}M$, the probability distribution of the solution starting from $x_0$.
An example 
is the gradient system induced from an isometric immersion 
$\alpha : M\to  \R ^m$ with $X(x): \R ^m\to T_xM$ 
defined to be the orthogonal projection for each $x\in M$. 
Another class of examples arise from symmetric space structures on $M$, 
see \cite{Elworthy-LeJan-Li-book}.

For our fixed $x_0$ in $M$ there is the solution map, or {\it It\^{o} map}, 
$$\ {\cal I}:C_0 \R ^m\to C_{x_0}M,$$ 
of (\ref{1}) defined by 
\[ 
{\cal I}(\omega )_t=x_t(\omega ),\qquad \omega \in C_0\R ^m, 
\] 
where $x_t$ is the solution starting at $x_0$. Thus ${\ {\cal I}%
_{*}({\Bbb P})=\mu _{x_0}}$. This It\^{o} map has an H-derivative in the 
sense of Malliavin calculus which is a continuous linear map from the Cameron-Martin space $H\equiv L_0^{2,1}\R^m$,
\[ 
{T_\omega {}}{\cal I}: H \to T_{\I(\omega )}C_{x_0}, 
\] 
for almost all $\omega \in C_0\R ^m$. Thus for $h\in  H $ 
and $0\le t\le T$ we have $T{\cal I}(h)_t\in T_{x_t}M$, a.s. 
 
\bigskip 
 
{\bf B. } Let ${\{\xi _t:0\le t\le T\}}$ denote the flow of (%
\ref{1}) so ${{\cal I}(\omega )_t=x_t(\omega )=\xi 
_t(x_0,\omega )}$. It can be taken to consist of random $C^\infty $ 
diffeomorphisms ${\xi _t:M\to M}$ with derivative maps $%
{T\xi _t:TM\to TM}$, so that ${T_{x_0}\xi _t\in 
\L  \left( T_{x_0}M;T_{x_t}M\right) }$.

Take ${h\in  H }$. Set ${v_t=T%
{\cal I}(h)}_t$. Bismut showed that $v_{\cdot }$ satisfies the 
covariant equation along the paths of ${\{x_t:0\le t\le T\}}$ 
\begin{equation} 
Dv_t=\nabla _{v_t}X\circ dB_t+X(x_t)\dot{h}_tdt  \label{2} 
\end{equation} 
with solution 
\begin{equation} 
v_t=T\xi _t\int_0^t(T\xi _s)^{-1}\left( X(x_s)\dot{h}_s\right) ds  \label{3} 
\end{equation}

\begin{lemma}  \label{le:decomposition-of-noise} \cite{Elworthy-Yor, Elworthy-LeJan-Li-book} There is a 
canonical decomposition of the noise $\{B_t:0\le t\le T\}$ given by 
\begin{equation} 
dB_t=\tilde{\parals_t}d\tilde{B}_t+\tilde{\parals_t}d\beta _t  \label{4} 
\end{equation} 
where 
 
\begin{enumerate} 
\item[(i)]  $\{\tilde{B}_t:0\le t\le T\}$ is a Brownian motion on the 
orthogonal complement $[\ker X(x_0)]^{\perp }$ of the kernel of $X(x_0)$ in $\R ^m$; 
 
\item[(ii)]  $\{\beta _t:0\le t\le T\}$ is a Brownian motion on $\ker X(x_0)$; 
 
\item[(iii)]  for each $t\ge 0$, 
$\tilde{\parals_t}:C_{x_0}M\to O(m)$ is a measurable map into the 
orthogonal group of $ \R ^m$ with $\tilde{\parals_t}(\sigma)[\ker X(x_0)]   =\ker X(\sigma _t)$ for $\mu _{x_0}$ almost all $\sigma \in C_{x_0}M$. 
\end{enumerate} 
\end{lemma} 
 
{\bf N.B.} We will regularly consider random variables on ${%
C_{x_0}M}$, such as ${\tilde{\parals_t}}$, to be 
 random variables on 
${C_0 \R ^m}$ by taking their composition with ${\cal I}$. 
For example the stochastic equation (\ref{4}) above is to be interpreted 
that way. Moreover let ${{\cal F}^{x_0}}$ be the $\sigma $-algebra on $C_0 \R ^m$ generated by ${\cal I}$ with $
\{\F_t^{x_0}, 0\le t\le T\}$ the  filtration generated 
by $(x_s:0\le s\le T)$. Then we can, and often will, 
consider $\F^{x_0}$-measurable functions as functions, 
defined up to equivalence, on $C_{x_0}M$. 
 
\bigskip 
 
Let $\F^{\beta}$ be the $\sigma $-algebra generated by $\{\beta _t: 0\le t\le T\}$, and $\F^{\tilde B}$ that generated by $\{ \tilde B_t: 0\le t\le T\}$. From Elworthy-Yor \cite{Elworthy-Yor}, 
Elworthy-LeJan-Li \cite{Elworthy-LeJan-Li-book} we know that 
 
\begin{itemize} 
\item[(a). ]  $\F^\beta $ and $\F^{\tilde{B}}$ are independent 
and 
 
\item[(b). ]  $\F^{\tilde{B}}$ $=\F^{x_0}$. 
 
\item[(c). ]  Equation (\ref{2}) can be written as the It\^{o} equation 
\begin{equation} 
Dv_t=\nabla _{v_t}X\left( \tilde{\parals_t}d\beta _t\right) -\frac 
12 {\Ric}^{\#}(v_t)dt+X(x_t)\dot{h}_t\,\,dt. 
\label{filter} 
\end{equation} 
where ${\Ric}_x^{\#}:T_xM\to T_xM$ corresponds to the Ricci 
curvature by  
$\langle {\Ric}^{\#}_x(u^1),u^2\rangle =\Ric(u^1,u^2)$ for $%
u^1,u^2$ in $T_xM$. 
\end{itemize} 
 We shall often write covariant derivatives such as $\nabla _vX$ as $\nabla 
X(v)$ so $\nabla X\left( v\right) \circ dB_t$ is just $\nabla _vX\circ dB_t$.

{\bf D.} We first show that $\wedge^q T_\omega \I$ take values in the 
exterior product space $\wedge^q T\C_{x_0}$ rather than just in
$\wedge^q_\epsilon T\C_{x_0}$.
Recall that a continuous linear map of  $H$ to a 
separable Banach space $E$ is ${\gamma}$-{\it radonifying} if it maps the  canonical Gaussian cylinder set measure of $H$ to a Borel measure on $E$. The {\it  2-summing norm} of an operator $A: E\to E$ is
$$\pi_2(A)=\sup_{\{x_n\}\subset E}
{\sum  \|A x_n\|^2
\over \sup_{\|u\|=1, u\in E^*}\sum \big( u (x_n)\big) ^2}.$$
 
\begin{lemma} \label{2B} 
For almost all $\omega \in C_0 \R ^m$ the map 
$$T_\omega {\cal I}\hskip 6pt:\hskip 6pt H  
  \longrightarrow T_{{\cal I(\omega )}}C_{x_0}$$ 
is $\gamma $-radonifying. Its operator norm $\|T{\cal I}\|$ is in $L^p(C_0%
 \R ^m)$ for $1\le p<\infty $ as is the 2-summing norm of its adjoint. 
\end{lemma}

\noindent{\it Proof.}
 Note that $\alpha: h\mapsto \int_0^{\cdot }(T\xi 
_s)^{-1}X(x_s)(\dot{h}_s)ds$ maps $H$ to $L_0^{2,1}(T_{x_0}M)$ and is continuous linear; almost surely. 
The inclusion $i: L_0^{2,1}(T_{x_0}M)\to C_0T_{x_0}M$ is $\gamma $-radonifying.  Write $T\I=T\xi_\cdot \circ i\circ \alpha$. Then the first result follows by composition properties of $\gamma $-radonifying maps and continuity of $T_{x_0} \xi 
_{\cdot }: C_0T_{x_0}M\to T_{x_\cdot(\omega )}C_{x_0}$. The 
p-th power integrability of the operator norms come from the corresponding  properties of $T\xi _t$ and $(T\xi _t)^{-1}$, e.g. see 
Kifer \cite {Kifer88}.  For the 2-summing norm apply Schwartz's duality theorem  \cite{Schwartz69} to see that the adjoint of the  $\gamma$-radonifying map $i$  is  $2$-summing with norm $\pi_2(i)$. 
 Then use the composition properties  of $2$-summing operators to estimate the $2$-summing norm
  $$\pi_2((T\I)^*)\le \|\alpha^*\|\; \pi_2(i^*)\;\|(T\xi_t)^*\|, \quad a.s.$$
Then apply  the integrability results again to see the norm is in $L^p$.
\hfill\rule{2mm}{2mm}

\begin{theorem} 
For almost all $\omega $ the map $\wedge ^qT_\omega I$ can be considered as 
a continuous linear map 
\[ 
\wedge ^q\left( T_\omega {\cal I}\right) :\wedge ^q\left( L_0^{2,1} \R %
^m\right) \longrightarrow \wedge ^q\left( T_{{}{\cal I}(\omega 
)}C_{x_0}\right) . 
\] 
Moreover the operator norms lie in $L^p(C_0 \R ^m)$ for $1\le p<\infty $. 
\end{theorem}

\noindent{\it Proof.}
This follows from Lemma \ref{2B} and results of Carmona \& 
Chevet \cite{Carmona-Chevet79} especially their Proposition 3.1 and Lemmas 3.1  a version of which is stated below as Lemma \ref{Carmona-Chevet}. Although they only deal with tensor products of two maps the lemma shows that the result holds for general $q$ by induction. 
\hfill\rule{2mm}{2mm}

\begin {lemma}[Carmona \& Chevet] 
\label{Carmona-Chevet}
Consider separable Hilbert spaces $H$ and $K$ and separable Banach spaces $E$ and $F$.  Let $T: H\to E$ be $\gamma$- radonifying and $S:K\to F$ bounded linear.
Then $S\otimes T:H\otimes K\to E\otimes _\pi F$ is a bounded linear map into the projective tensor product. Moreover $$ \|S\otimes T\|_{\mathbb{L}(H\otimes K;E\otimes_\pi F)}\leqslant \pi_2(T^*)\|S\|$$
where $\pi_2(T^*)$ denotes the 2-summing norm of the adjoint of $T$.
\end {lemma}

\medskip

The conditional expectations of these operators can be defined as in 
Elworthy-Yor  \cite {Elworthy-Yor}, Elworthy-LeJan-Li 
\cite{Elworthy-LeJan-Li-book}, to give bounded linear maps, 
defined almost surely, 
$$\overline{\wedge ^q\left( T {\cal I}\right)} (\omega) :\wedge ^q 
 H  \longrightarrow \wedge ^q 
 \left( T_{{ \cal I}(\omega )} C_{x_0}\right).  $$ 
For example 
$$\overline{\wedge ^q\left( T{\cal I}\right) } (\omega) 
:= \E \left\{ \wedge 
^q\left( T_\omega {\cal I}\right) \,\, \left|\;  \F^{x_0}\right. \right\} (\omega ) $$ 
is given by 
$$\overline{\wedge ^q\left( T{\cal I}\right) } (\omega)(h)_t\,=\,\,\left( 
\wedge ^q \parals_t\right) 
 \,\, \E \,\left\{ \,\wedge ^q\left( \parals_t^{-1}\right) 
\,\wedge ^q\left( T{\cal I}_t\right (h)\,\,\left|\,\F^{x_0} 
\right. \right\}(\omega) .  $$ 
 
For $\mu _{x_0}$ almost all $\sigma \in C_{x_0}M$ we have also 
$$\left(\overline{\wedge ^q\left( T{\cal I}\right) }\right) 
_\sigma:\wedge ^q H \longrightarrow \wedge ^q\left( T_\sigma C_{x_0}\right)$$ 
given by 
$$\left(\overline{\wedge ^q\left( T{\cal I}\right) }\right)_\sigma (h) 
\,\,:=\,\, \E %
\,\left\{ \,\wedge ^q\left( T_{\cdot }{\cal I}\right) (h)\,\,\left| \,\,{}%
{\cal I}=\sigma \right. \,\right\}.$$ 
Note the inequalities 
\begin{eqnarray*} 
\left\|\overline{\wedge ^q\left( T{\cal I}\right) } (\omega)(h)\right\| 
 &\le & \E  \left\{ 
 \left|\wedge ^q\left( T_{\cdot }{\cal I}\right)(h)\right| 
\,\,\left\vert \,{\cal F}^{x_0}\right. \right\} (\omega )\qquad a.s. \\ 
&\le & \E \left\| \wedge ^q\left( T_\omega {}{\cal I}\right) \right\| 
\; |h| \qquad a.s. 
\end{eqnarray*} 
which give $L^p$ bounds for operator norms of ${\overline{%
\wedge ^q\left( T{\cal I}\right) }}$, $q=1,2,\dots $. 
 
\bigskip 
 
{\bf E. Definition of $\HH_\sigma ^q$, $H$-$q$-vector fields and $H$-$q$-forms.}  We can now define $  \HH_\sigma ^q$, for almost all
$\sigma \in C_{x_0}M$, to  be the image of  
$\overline{\wedge ^q\big( T{\cal I}\big) }_\sigma $ 
 in $\wedge ^qT_\sigma C_{x_0}$ together with the 
inner product induced by the linear bijection 
 
\[ 
{\overline{\wedge ^qT\I}_\sigma }_{\left| _{[\ker\overline{\wedge 
^qT\I}_\sigma ]^{\perp }}\right. }: \quad [\ker\overline{%
\wedge ^q\left( T\I\right) }_\sigma ]^{\perp }\to \HH_\sigma ^q. 
\] 
Thus the ${\HH _\sigma ^q}$ are Hilbert spaces with natural  continuous linear inclusions ${\iota _\sigma }$, say, into the $T_\sigma C_{x_0}$.

 Denote by $\HH ^q=\cup_\sigma \HH _\sigma^q$ the ``vector bundle over 
$C_{x_0}M$" with fibres $\HH _\sigma^q$, and $(\HH ^q)^*$ 
the corresponding dual ``bundle". Set $\cal H =\HH ^1$. Since these are only almost sure defined it is not strictly speaking correct to consider them as bundles over $C_{x_0}M$ though some vector bundle structure is given to $\cal H$ in \cite {Elworthy-Li-Ito-map} see also \S\ref{se-torsion} below.
The space of $L^2$ sections of $\HH ^q$ 
and ${\HH ^q}^*$ are denoted by $L^2\Gamma{\HH ^q}$ and 
$L^2\Gamma{\HH ^q}^*$. Sections of $(\HH ^q)^*$  or of $(\HH ^q)$ will be called  $H$-$q$-forms (or admissible $q$-forms), or $H$-$q$-vector fields respectively. Note that any q-form on $C_{x_0}M$ restricts to give an $H$-$q$-form.

\section{Characterization of $\HH _\cdot^1$ and $\HH _\cdot^2$} 
 \label{se-characterization} 
 
{\bf A.} `Damped parallel translations' $W_t^{(q)}$ will play an essential 
role. For a $q$-vector ${\ v\in \wedge ^qT_{x_0}M}$, define $W_t^{(q)}(V)\in \wedge ^qT_{x_t}M$ to be the random $q$%
-vector satisfying 
\begin{equation} 
{\frac D{dt}}W_t^{(q)}(V)=-{\frac 12}{\cal R}^q\,W_t^{(q)}(V),\hskip %
30pt0\le t\le T  \label{damped} 
\end{equation} 
where ${\cal R}^q\in \Hom(\wedge ^qTM;\wedge ^qTM)$ is the 
Weitzenbock curvature term deinfed by ${\cal R}=\Delta-\trace \nabla^2$, see e.g.  Airault \cite{Airault76}, 
Elworthy\cite{ELflour}, Ikeda-Watanabe\cite{Ikeda-Watanabe}, 
Elworthy-LeJan-Li 
\cite{Elworthy-LeJan-Li-book}, Elworthy-Li-Rosenberg \cite{Elworthy-Li-Rosenberg}, Malliavin \cite{Malliavin-book}. Here (\ref{damped}) is a covariant equation along the paths of our solution $\{x_t:0\le t\le T\}$ to (\ref {1}). 
 
\bigskip 
 
For $q=1$ write $W_t$ for $W_t^{(1)}$. Then ${W_t:T_{x_0}M\to 
T_{x_t}M}$ is the Dohrn-Guerra translation given by 
$${\frac D{dt}}W_t(V)=-{\frac 12}{Ric}_{x_t}^{\#}\left( W_t(V)\right) ,\qquad    0\le t\le T.$$
Write 
\[ 
{{\frac{{\rm I\!D}}{dt}}=W_t{\frac D{dt}}W_t^{-1}} 
\] 
acting on suitably regular vector fields ${\{v_t:0\le t\le T\}}$ 
along the paths of ${\{x_t:0\le t\le T\}}$. Then 
$${{\frac{{\rm I\!D}}{dt}}={\frac D{dt}}+{\frac 12}Ric^{\#}},$$ 
c.f. Fang, formula (1,3), in Fang \cite{Fang98} and Norris \cite{Norris-twisted}. 
 
\begin{definition}\label{L2T}
For almost all paths $\omega$, define the $L^2$ tangent space $ L^2T_{\sigma_{\cdot }}C_{x_0}$  to  consist of those paths $u:[0,T]\to TM$  over $\sigma $ with 
$$\parals_{\cdot }^{-1}u_{\cdot }\in L^2\left( 
[0,T];T_{x_0}M\right) $$ together with its natural Hilbert space structure. 
\end{definition}
It was shown in \cite{Elworthy-LeJan-Li-CR}, see also  \cite{Elworthy-Li-Hodge-1}, \cite{Elworthy-LeJan-Li-book} that 
\begin{equation} 
\overline{T{\cal I}{}}_t(h)={\cal W}_t\left( X(x_{\cdot })\dot{h}_{\cdot 
}\right)  \label{Ito-map} 
\end{equation} 
where 
\[ 
{\cal W}:L^2T_{x_{\cdot }}C_{x_0}\to T_{x_{\cdot }}C_{x_0} 
\] 
is defined by 
\begin{equation} 
\left( {\cal W}(u)\right) _t=W_t\int_0^t\left( W_r\right) ^{-1}\left( 
u_r\right) dr.  \label{w} 
\end{equation} 
Note that
\begin{equation} 
{\frac{{\rm I\!D}}{dt}}\left( {\cal W}{}(u)\right) _t=u_t, \qquad
u\in L^2T_{x_{\cdot }}C_{x_0}.  \label{Wei-map-1} 
\end{equation} 

Thus, as shown in \cite{Elworthy-LeJan-Li-book}, \cite{Elworthy-Li-Hodge-1}, 
\begin{equation}\label{eq-H1}
\HH _\sigma ^1=\left\{ v\in T_\sigma C_{x_0}:\parals_{\cdot }^{-1}v_{\cdot 
}\in L_0^{2,1}\left( T_{x_0}M\right) \right\} 
 \end{equation}
with inner product 
\begin{equation} 
\langle v^1,v^2\rangle _{\HH ^1}=\int_0^T\left\langle {\frac{{\rm I\!D}}{ds}}v_s^1,{\frac{{\rm I\!D}}{ds}}v_s^2\right\rangle ds 
\end{equation} 
so that ${\D\over d \cdot}:\HH _\sigma ^1\to L^2T_\sigma 
C_{x_0}$ is an isometric isomorphism with inverse $\cal W$ for almost all $\sigma \in C_{x_0}M$. Thus it agrees as a Hilbertable space with the usual Bismut tangent space, though the inner product is not the one originally used. 
Using the same notation, by \S\ref{se-exterior}D  we note that a vector 
 $u$ of  $\wedge^2T_\sigma C_{x_0}M$ is in  $\wedge^2\HH _\sigma$ 
if and only if there exists $\underline {k}\in \wedge^2 L^2T_\sigma C_{x_0}M$ 
so that 
$$u_{s,t}=(\wedge^2 {\cal W})_{s,t}\underline k, $$
or written in full,
\begin{equation}u_{s,t}=  \left(W_t\int_0^t(W_{r_1})^{-1}(-) dr_1 
 \wedge W_t\int_0^t(W_{r_2})^{-1}(-) dr_2 
\right)\underline k_{r_1,r_2}. \end{equation} 
If so $\underline k_{s,t} 
 ={\D\over \partial s}\otimes{\D\over \partial t}u$ 
or equally 
$\underline k=\wedge^2 {\D\over d\cdot}u$. 
 
\bigskip

{\bf B. } 
 More generally let
  ${L^2\left( \wedge ^qTM\right) _\sigma }$ and $%
{C_0\left( \wedge ^qTM\right) _\sigma }$ denote respectively  
the spaces of $L^2$ and continuous paths vanishing at $0$, ${u:[0,T]\to \wedge 
^qTM}$ over $\sigma $. 
 Define $${\cal W}^{(q)}:L^2\left( \wedge ^qTM\right) _\sigma
  \to C_0\left( \wedge ^qTM\right) _\sigma $$ 
 by 
\begin{eqnarray} 
\left( {\cal W}^{(q)}(V_{\cdot })\right) _t &=&W_t^{(q)}\int_0^t\left( 
W_r^{(q)}\right) ^{-1}(V_r)dr  \label{Wei-q-formula} \\ 
&=&\int_0^t{W^{(q)}}^r_t\left( V_r\right) dr, 
\end{eqnarray} 
where 
$${W^{(q)}}^s_t=W_t^{(q)}\left( W_s^{(q)}\right) ^{-1}$$ 
is the solution to 
\begin{equation}
\label{Weitzenbock-q}
{\frac D{dt}}{W^{(q)}}^s_t(V) 
=-{\frac 12}{\cal R}^q\left({W^{(q)}}^s_t(V)\right),\qquad   s\le t\in[0,T]
\end{equation}
with ${W^{(q)}}^s_s=\id:\wedge ^qT_{\sigma _s}M\rightarrow \wedge ^qT_{\sigma 
_s}M$. Write $W^s_t$ for ${W^{(1)}}^s_t$ and observe that ${\cal W}^{(1)}={\cal W}$.  For simplicity we shall 
 write ${\cal W}_t^{(q)}(V_{\cdot })$ for 
 $\left( {\cal W}^{(q)}(V_{\cdot })\right)_t$. 
 
Set 
\begin{equation} 
\label{big-D} 
\frac{\D ^{(q)}}{dt}=\left( \frac D{dt}\right) +\frac 12{\cal R}^q, 
\end{equation} 
acting on $q$-vectors on $M$ along a sample path $\sigma $. Then as for 
 $q=1$, and for $W_t^{(q)}$ defined by (\ref{Weitzenbock-q}): 
$$\frac{\D^{(q)}}{dt}V_{t,\dots ,t}=W_t^{(q)}\frac 
d{dt}(W_t^{(q)})^{-1}V_{t,\dots ,t} $$
and the inverse of $\frac{\D ^{(q)}}{d\cdot}$ is 
$$\left(\frac{\D ^{(q)}}{d\cdot}\right)^{-1}=W_\cdot^{(q)}\int_0^\cdot W_r^{(q)}(\ev_r-)dr={\cal W}^q$$
where $\ev_r$, generically, denotes the evaluation operator at $r$.
Furthermore let ${\cal R}:\wedge ^2TM\rightarrow \wedge ^2TM$ be the  curvature operator. Then the second Weitzenbock curvature ${\mathcal R}^2$ is given by
$${\cal R}^2=d\wedge ^2\left( \Ric^{\#}\right) -2{\cal R}.  $$ 
Therefore using (\ref{section-1-1}), for $V\in \wedge^2T_\sigma C_{x_0}M$, 
\begin{equation} 
\frac{\D ^{(2)}}{dt}V_{t,t} 
 =\left( \left( (d\wedge^2 )\big(\frac{\D }{d\cdot }\big)\right) V\right)_{t,t} 
 -{\cal R}\left( V_{t,t}\right),
\label{relate-R2-R} 
\end{equation} 
whenever all the terms involved make sense.
In the above we have identified ${\frac D{dt}}V_{t,t}$ with 
 $\left( d\wedge ^2{\frac D{d\cdot }}\right) (V)_{t,t}$ where the first 
 refers to covariant  differentiation of the 2-vector field 
 $\{V_{t,t}:0\le t\le T\}$ along $\sigma $ obtained from the element $V$ in 
 $\wedge ^2T_\sigma C_{x_0}$. 

\bigskip 
 
{\bf C.} In this section we shall discuss a system of equations related to 
the conditional expectation of the It\^{o} map. First note that  the 
curvature operator ${\cal R}$ on the manifold $M$ induces a linear map 
$Q_\sigma$ on  $\wedge ^2_\epsilon T_\sigma C_{x_0}$ given  by 
\begin{equation}\label{Q} 
Q_\sigma (G)_{s,t}=\left( \1\otimes W_t^s\right) 
 W_s^{(2)}\int_0^s\left( W_r^{(2)}\right) ^{-1}\left({\cal 
R}_{\sigma _r}\left( G_{r,r}\right) \right) dr,\,\,\,\,\,\,s\le t. 
\end{equation} 
Equivalently, 
\[ 
Q_\sigma (G)_{s,t}=
\left(W_s\otimes W_t\right)
\left(\wedge^2 (W_\cdot^{-1}) W_\cdot^{(2)}\int_0^\cdot\left( W_r^{(2)}\right) ^{-1} 
\left({\cal R}_{\sigma _r}\left( G_{r,r}\right) \right) dr\right)_{min(s,t)}. 
\] 
Clearly 
\begin{equation}\left\{
\begin{array}{ccl}
(\1\otimes {\D\over dt})Q(G)_{s,t}&=&0, \hskip 18pt s<t,\\
{\D^{(2)}\over ds}Q(G)_{s,s}&=&{\mathcal R}(G_{s,s}).
\end{array}\right. .
\end{equation}
The second equation is equivalent to
 $\left(d\wedge^2{\D\over ds}\right)Q(G)_{s,s}={\mathcal R}((I+Q)G)_{s,s}$. 
Define $j_G: [0,T]\to T_{x_0}M\otimes T_{x_0}M$ by 
\begin{equation} \label{eq-jay-def}
j_G\left( s\right) =\left( W_s^{-1}\otimes W_s^{-1}\right) 
W_s^{(2)}\int_0^s\left( W_r^{(2)}\right) ^{-1}\left( {\cal R}_{\sigma _r} 
\left( G_{r,r}\right) \right) dr. 
\end{equation} 
Then  $j_G$ is $C^1$ and\begin{equation}\label {eq-jay}
\left( W_s^{-1}\otimes W_t^{-1} \right)Q_\sigma \left( G\right)_{s,t}
=j_G\left( s\wedge t\right).
\end{equation}

If we set
\[ 
D\left( \wedge^2 T_\sigma C_{x_0}\right) =\left\{ \begin{array}{l} 
  u\in \wedge ^2_\epsilon T_\sigma C_{x_0} \hbox{ such that }\\ 
(1)\; \hbox{for each } 0\le s<T, 
   ,\,\,t\mapsto \left( \parals_s^{-1}\otimes \parals_t^{-1}\right) 
 u_{s,t} \hbox{ is }\\ 
   \hbox{  absolutely  continuous on } $(s,T]$; \\ 
(2)\; r\mapsto \wedge ^2(\parals_r^{-1})u_{r,r}\hbox{ is absolutely 
     continuous on } $[0,T]$ \\ 
\end{array} 
\right\} 
\] 
then $Q(G)$ clearly lies in $D( \wedge^2 T_\sigma C_{x_0})$.
There is another linear map $\pathR$  on  $\wedge ^2_\epsilon T_\sigma C_{x_0}$
 defined by:
\begin{equation}\label{eq:R}
\pathR (Z)_{s,t}=\left( W_s\otimes W_t\right) \int_0^s\left( \wedge 
^2W_r{}^{-1}\right) \left( {\cal R}_{\sigma _r}\left( Z_{r,r}\right) 
\right) dr, \hskip 10pt  s\le t,
\end{equation} 
which also sends $\wedge^2_\epsilon T_\sigma C_{x_0}M$ to 
 $D( \wedge^2 T_\sigma C_{x_0}M)$. Furthermore, from equation ( \ref{relate-R2-R})
\begin{equation} \left\{
\begin{array}{ccl}
(\1\otimes {\D\over dt})\pathR(Z)_{s,t}&=&0, \hskip 18pt s<t,\\
{\D^{(2)}\over ds}\pathR(Z)_{s,s}&=&{\cal R}_{\sigma_s}\left(Z_{s,s}-\pathR(Z)_{s,s}\right).
\end{array}\right. 
\end{equation}
In fact $\1+Q$ and $\1-\pathR$ are inverse of each other as described in 
the following lemma. It will be shown later, \S\ref{se-torsion}.2, that $\pathR$ restricted to 
$\wedge^2 {\mathcal H}^1$ is the curvature operator of the damped Markovian
 connection on  ${\mathcal H}^1$ which is induced by the map 
${\D\over d\cdot}$  from the pointwise connection on  the $L^2$ tangent bundle $L^2TC-{x_0}$. 
 
\begin{lemma} 
\label{lemma-chara-1} 
\begin{enumerate} 
\item[(i)]  Given $G\in D\left( \wedge^2 T_\sigma C_{x_0}\right) $, there 
 is a unique  solution $Z\in D\left( \wedge^2 T_\sigma C_{x_0}\right)$, 
 to the following equations 
 
\begin{equation} 
\left\{ 
\begin{array}{cll} 
\left( \1\otimes \frac{\D }{dt}\right) Z_{s,t}&=&\left( \1\otimes \frac{%
\D }{dt}\right) G_{s,t}, \hskip 15pt s< t, \\ 
\frac{\D ^{(2)}}{ds}Z_{s,s}&=&\left(\left( \left( d\wedge ^2\right) \left(\frac{{\Bbb %
D}}{d\cdot }\right)\right)G \right) _{s,s}\\ 
Z_{0,0}&=&G_{0,0}. 
\end{array} \right. 
\label{characteristic-eq-1} 
\end{equation} 
The solution is: 
\[ 
Z_{s,t}=G_{s,t}+Q_\sigma(G)_{s,t}. 
\] 
Conversely for each $Z\in   D\left( \wedge^2 T_\sigma C_{x_0}\right) $ the
unique solution to (\ref{characteristic-eq-1} ) is given by:
\begin{equation} \label {eq-G-Z}
G=Z-{\pathR}(Z).\end{equation} 
 
\item[(ii)] As operators on $\wedge^2_\epsilon T_{\sigma}C_{x_0}M$  both $Q$ and $\pathR$ are compact and  $\1+Q$ and $\1-\pathR$ are mutual inverses. In particular for all $v$ in
 $\wedge^2_\epsilon T_{\sigma}C_{x_0}M$,
$$Q(v)=\pathR(v+Q(v)),$$
$$Q(\1+Q)^{-1}v=\pathR(v).$$
\item[(iii)]  The following holds on $ D\left( \wedge^2 T_\sigma C_{x_0}\right)$:
\begin{equation} 
(\wedge^2 W^{-1}Z)_{s,t}- (\wedge^2 W^{-1}Z)_{s\wedge t} 
=(\wedge^2 W^{-1}G)_{s,t}-(\wedge^2 W^{-1}G)_{s\wedge t},  \label{section-3-1} 
\end{equation} 
which is equivalent to, for $r\le s\le t$, 
\[ 
Z_{r,t}-\left( \1\otimes W_t^s\right) Z_{r,s}=G_{r,t}-\left( 
\1\otimes W_t^s\right) Z_{r,s}. 
\]
\end{enumerate} 
\end{lemma} 
 
\noindent 
{\it Proof}. Given $G\in  D\left( \wedge^2 T_\sigma C_{x_0}\right) $,
 $Z=(\1+Q_\sigma)(G)$ certainly solves  (\ref {characteristic-eq-1}). 
For uniqueness let $Z$ be any solution in $D( \wedge^2 T_\sigma C_{x_0})$. 
 Solve the first equation in (\ref{characteristic-eq-1}) to get 
\begin{equation} 
Z_{s,t}=G_{s,t}+\left( W_s\otimes W_t\right) \left( \tilde{j}(s)\right) 
,\,\,\,s\le t  \label{characteristic-eq-3} 
\end{equation} 
some $\tilde{j}(s)\in \wedge ^2T_{\sigma _0}M$. Then 
\begin{equation} 
Z_{s,s}=G_{s,s}+\left( W_s\otimes W_s\right) \left( \tilde{j}(s)\right) . 
\label{section-3-2} 
\end{equation} 
In particular $\left( W_s\otimes W_s\right) \left( \tilde{j}(s)\right) $ 
 is absolutely continuous in $s$. Substitute the above equation 
(\ref{section-3-2}) 
into (\ref{characteristic-eq-1}) and use (\ref{relate-R2-R}) to see 
$$\frac{\D ^{(2)}}{ds}\left( W_s\otimes W_s\right) \left( \tilde{j}%
(s)\right) ={\cal R}\left( G_{s,s}\right), $$ 
giving 
$$\left( W_s\otimes W_s\right) \left( \tilde{j}(s)\right) 
=W_s^{(2)}\int_0^s\left( W_r^{(2)}\right) ^{-1}\left({\cal 
R}_{\sigma _r} \left( G_{r,r}\right) \right) dr.  $$
Thus $ \tilde{j}(s) =j_G(s)$ and uniqueness holds by formula (\ref{eq-jay}).

 Similarly given $Z\in  D\left( \wedge^2 T_\sigma C_{x_0}\right) $,
$(\1-\pathR)(Z)$ is seen to satisfy (\ref{characteristic-eq-1}) given $Z\in  D\left( \wedge^2 T_\sigma C_{x_0}\right) $.

Now using the isometry between $\otimes^2_\epsilon C_0T_{x_0}M$ and $C_0\left( [0,T]^2;\otimes ^2 T_{x_0}M\right)$ and the Arz\`ela-Ascoli theorem applied to $(s,t)\mapsto j(s,t)$ for a bounded set of $G$, we see that $Q:\wedge^2_\epsilon T_{\sigma}C_{x_0}M\to\wedge^2_\epsilon T_{\sigma}C_{x_0}M$ is compact. Therefore $\1+Q$ has closed range. Since we have just seen that its range contains all $Z$ in the dense subspace $D\left( \wedge^2 T_\sigma C_{x_0}\right) $ it is surjective and so an isomorphism. By equation (\ref{eq-G-Z}) its inverse is $\1-{\pathR}$ and so $\pathR$ is compact. The rest of parts (i) and (ii) follow directly.

Part (iii) follows from (\ref {characteristic-eq-3}) and (\ref{section-3-2}).\hfill \rule{2mm}{2mm} 
 
See \S\ref{se-tensor} below for a more detailed examination of $Q(V)$. 
\bigskip

{\bf D.} 
The following theorem gives alternative descriptions of the space 
$\HH ^2_\sigma$. 
\begin{theorem} 
\label{Theorem-3.2} 
For any $h^1, h^2\in L_0^{2,1}\R^m$, set $\underline h=h^1\wedge h^2$. Then 
\begin{equation}\label{wedge-2} 
\overline{\wedge^2T\I}({\underline h})=(\1+Q)\wedge^2\overline 
{T\I}({\underline h}). 
\end{equation} 
In particular the space $\HH _\sigma ^2=\{\overline{\wedge ^2T%
{\cal I}}_{\sigma}(h),h\in \wedge ^2 H \}$ can be characterised by 
either of the following: 
\begin{enumerate} 
\item[(i)] 
\[ 
\HH _\sigma ^2=\left\{     \begin{array}{l} 
u\in D(\wedge^2T_\sigma C_{x_0}),\hbox{such that there exists }G\in \HH %
_\sigma ^1\wedge \HH _\sigma ^1,\hbox{with} \\ 
\left( \left( \1\otimes \frac{\D }{d\cdot }\right) u\right) 
_{s,t}=\left( \left( \1\otimes \frac{\D }{d\cdot }\right) 
G\right) 
_{s,t},\,\,\,\,\,\,\,s< t, \\ 
\hbox{and }\; \frac{\D ^{(2)}}{ds}u_{s,s} 
=\left( \left( \left( d\wedge ^2\right) 
\frac{\D }{d\cdot }\right) G\right) _{s,s}\,\,\,\,\,\,\,0\leqslant s\leqslant T
\end{array} 
\right\} . 
\] 
\item[(ii)] 
\[ 
\HH _\sigma ^2=\left\{ u\in \wedge ^2_\varepsilon T_\sigma C_{x_0},\hbox{such that } 
u=v+Q_\sigma(v),\hbox{ some }v\in \HH _\sigma ^1\wedge 
 \HH _\sigma ^1\right\},  \] 
and for $v_1, v_2\in \wedge^2\HH ^1_\sigma$,  by definition, 
\begin{equation} 
\langle v^1+Q_\sigma(v^1), v^2+Q_\sigma(v^2)\rangle_{\HH ^2_\sigma} 
=\langle v^1, v^2\rangle_{\wedge^2\HH ^1_\sigma}. 
\end{equation} 

\item[(iii)] $u\in \HH^2$ if and only if $u-\pathR(u)\in \wedge^2\HH^1$. If
so $$\| u\|_{\HH^2}=\| u-\pathR(u)\|_{\wedge^2 \HH^1}.$$
\end{enumerate} 
In particular $\HH ^2_\sigma$ depends on the Riemannian structure of $M$ 
but \underline{not} the choice of stochastic differential equation 
 (\ref{1}) provided its 
LeJan-Watanabe connection in the sense of Elworthy-LeJan-Li 
\cite{Elworthy-LeJan-Li-Tani} is the Levi-Civita connection. 
\end{theorem} 
 
\noindent 
{\it Proof. } For $h^1\wedge h^2\in L_0^{2,1}\left( R^m\right) $, write $%
V^1\wedge V^2=\left( \wedge ^2T{\cal I}\right) \left( h^1\wedge h^2\right)$. 
 Then applying It\^{o}'s formula in $t$ for $0\le s< t\le T$ with $D_t$ 
referring to covariant stochastic differentiation in $t$, 
\begin{eqnarray*} 
&&\left( \1\otimes D_t\right) \left( V^1\wedge V^2\right) _{s,t} \\ 
&=&\frac 12V_s^1\otimes \left( \nabla X\left( V_t^2\right) \circ 
dB_t+X\left( x_t\right) \left( \dot{h}_t^2\right) dt\right) \\ 
&&-\frac 12V_s^2\otimes \left( \nabla X\left( V_t^1\right) \circ 
dB_t+X\left( x_t\right) \left( \dot{h}_t^1\right) dt\right) \\ 
&=&\frac 12V_s^1\otimes \left( \nabla X\left( V_t^2\right) \parals_td\beta 
_t-\frac 12Ric^{\#}\left( V_t^2\right) dt+X\left( x_t\right) \left( \dot{h}%
_t^2\right) dt\right) \\ 
&&-\frac 12V_s^2\otimes \left( \nabla X\left( V^1_t\right) \parals_td\beta _t-\frac 
12Ric^{\#}\left( V_t^1\right) dt+X\left( x_t\right) \left( \dot{h}%
_t^1\right) dt\right) \\ 
&=&\left( \1\otimes \nabla X(-)\parals_td\beta _t\right) \left( 
V^1\wedge V^2\right) _{s,t}-\left( \1\otimes \frac 
12Ric^{\#}\left( -\right) \right) 
\left( V^1\wedge V^2\right) _{s,t}dt \\ 
&&+\frac 12\left( V_s^1\otimes X\left( x_t\right) \left( \dot{h}_t^2\right) 
-V_s^2\otimes X\left( x_t\right) \left( \dot{h}_t^1\right) \right) dt. 
\end{eqnarray*} 
Write $\overline{V^1\wedge V^2}$ for the conditional expectation of $%
V^1\wedge V^2$ with respect to ${\cal F}^{x_0}$, and similarly let $%
\overline{V}^i$ stand for the conditional expectation of $V^i$ with respect 
to ${\cal F}^{x_0}$. Then by (\ref{filter}), following Elworthy-Yor 
 \cite{Elworthy-Yor}, and (\ref{Wei-map-1}) 
 
\begin{eqnarray*} 
\left( \1\otimes D_t\right) \left( \overline{V^1\wedge V^2}\right) 
_{s,t} &=&-\left( \1\otimes \frac 12Ric^{\#}\left( -\right) 
\right) \left(\overline{ V^1\wedge 
V^2}\right) _{s,t}dt \\ 
&&+\frac 12\left( \overline{V}_s^1\otimes X\left( \dot h_t^2\right) 
  -\overline{V}_s^2\otimes X\left( \dot h_t^1\right) \right) dt. 
\end{eqnarray*} 
This is equivalent to 
\begin{eqnarray*} 
\left( \1\otimes \frac{\D }{d\cdot }\right) \left( 
\overline{V^1\wedge V^2}\right) _{s,t} 
  &=&\frac 12\left( \overline{V}_s^1\otimes X\left( \dot 
h_t^2\right) -\overline{V}_s^2\otimes X\left(\dot  h_t^1\right) \right) \\ 
&=&\left( \1\otimes \frac{\D }{d\cdot }\right) \left( \overline{V}%
^1\wedge \overline{V}^2\right) _{s,t}. 
\end{eqnarray*} 
On the other hand It\^{o}'s formula, applied to the $2$-vector field $%
\left\{V_t^1\wedge V_t^2,\,\,\,\,0\le t\le T\right\} $ along $\sigma $ in $M$, gives 
 
\begin{eqnarray*} 
D_t\left( V_t^1\wedge V_t^2\right) &=&V_t^1\wedge \left( \nabla X\left( 
V_t^2\right) \circ dB_t+X\left( x_t\right) \left( \dot{h}_t^2\right) 
dt\right) \\ 
&&+\left( \nabla X\left( V^1_t\right) \circ dB_t+X\left( x_t\right) \left( \dot{h}_t^1\right) dt\right) \wedge V_t^2 \\ 
\end{eqnarray*} 
Change to Ito differentials and decompose the noise recalling that  $\nabla X$ vanishes on $[\ker X]^{\perp }$: 
\begin{eqnarray*} 
D_t\left( V_t^1\wedge V_t^2\right) &=&V_t^1\wedge \left( \nabla X\left( 
V_t^2\right) \left( \parals_td\beta _t\right) -\frac 12Ric^{\#}\left( 
V_t^2\right) dt+X\left( x_t\right) \left( \dot{h}_t^2\right) dt\right) \\ 
&&+\left( \nabla X\left( V_t^1\right) \left( \parals_td\beta _t\right) -\frac 
12Ric^{\#}\left( V_t^1\right) dt+X\left( x_t\right) \left( \dot{h}%
_t^1\right) dt\right) \wedge V_t^2 \\ 
&&+\frac 12\sum_{i=1}^m\left( \nabla X^i\wedge \nabla X^i\right) \left( 
V_t^1\wedge V_t^2\right) \\ 
&=&\left( d\wedge ^2\left( \nabla X(-)\left( \parals_td\beta _t\right) \right) 
\right) \left( V_t^1\wedge V_t^2\right) \\ 
&&-\left( d\wedge ^2\left( \frac 12Ric^{\#}\left( -\right) \right) \right) 
\left( V_t^1\wedge V_t^2\right) dt \\ 
&&+ \sum_{i=1}^m\left( \nabla X^i\wedge \nabla X\right) ^i\left( 
V_t^1\wedge V_t^2\right) \\ 
&&+\left( V_t^1\wedge X\left( x_t\right) \left( \dot{h}_t^2\right) 
+X\left( x_t\right) \left( \dot{h}_t^1\right)\wedge V_t^2 \right) dt. 
\end{eqnarray*} 
But 
\begin{equation} 
-\left( d\wedge ^2\left( \frac 12Ric^{\#}\left( -\right) \right) \right) 
+
\sum_{i=1}^m\nabla X^i\wedge \nabla X^i=-\frac 12{\cal R}^2, 
\end{equation} 
as in Elworthy\cite{ELflow} for gradient systems, see also Elworthy-LeJan-Li 
\cite {Elworthy-LeJan-Li-book} for the general situation. Again use the 
 technique of Elworthy-Yor, \cite{Elworthy-Yor}, taking conditional 
 expectations to get: 
 
\[ 
\frac D{dt}\overline{V_t^1\wedge V_t^2}=-\frac 12{\cal R}^2\left( \overline{%
V_t^1\wedge V_t^2}\right) +\overline{V}_t^1\wedge X(x_t)\left( \dot h_t^2\right) +%
X(x_t)\left( \dot h_t^1\right)\wedge \overline{V}_t^2. 
\] 
Thus 
\begin{eqnarray*} 
\frac{\D ^{(2)}}{dt}\overline{V_t^1\wedge V_t^2} &=&\overline{V}%
_t^1\wedge X(x_t)\left( \dot h_t^2\right) 
+X(x_t)\left( \dot h_t^1\right)\wedge \overline{V}_t^2 
\\ 
&=&\left( d\wedge ^2\left( \frac{\D }{dt}\right) \right) \left( 
\overline{V}^1\wedge \overline{V}^2\right) _{t,t}. 
\end{eqnarray*} 
We have shown given $u=\overline{\wedge ^2T{\cal I}}\left( h^1\wedge 
h^2\right) $, it is related to 
 $\overline{T{\cal I}}\left( h^1\right) \wedge 
\overline{T{\cal I}}\left( h^2\right) $ by equation (\ref 
{characteristic-eq-1}). Solve the equation to obtain 
$$\overline{\wedge^2T\I}({\underline h})_{s,t} =\wedge^2\overline 
{T\I}({\underline h})_{s,t} +\left( \1\otimes W_t^s\right) 
W_s^{(2)}\int_0^s \left( W_r^{(2)}\right) ^{-1} \left( {\cal 
R}\left(\wedge^2 \overline{T\I}({\underline h})_{r,r}\right) 
   \right)dr,$$ 
 that is,  the desired identity (\ref{wedge-2}). On the other hand, 
 given $u$ satisfying (\ref{characteristic-eq-1}) for 
 $G=\wedge ^2\overline{T{\cal I}} (\underline h) $, 
$\underline h\in \wedge ^2L_0^{2,1}(\R^m)$, then 
 $u=\overline{\wedge ^2T{\cal I}}(\underline h) $ 
 by uniqueness of the solution. This proves the first equivalence. 
The second equivalence follows from Lemma \ref{lemma-chara-1}. Part (iii)
follows straightforwardly from the previous lemma.

\section{H-one-forms: exterior differentiation and Hodge decomposition}
\label{se-1-form}
{\bf A. Differentiation of functions.} For scalar analysis in our context and with this notation, we refer to \cite{Elworthy-Li-Ito-map} or for the basic facts to \cite {Elworthy-LeJan-Li-book}. As emphasised in \cite {Elworthy-Li-Ito-map} it is necessary to fix an initial domain,  $\Dom(d_\HH )\subset L^2(C_{x_0}M; \R)$ for  the H-derivative operator $d_\HH $. We shall choose this to be a subspace which contains the smooth cylindrical functions and consists of  $BC^2$ functions in the Fr\'echet sense, using the natural Finsler structure of $C_{x_0}M$, see \cite{Elworthy-Ma}).  For example the space of all smooth cylindrical functions. (We will require two derivatives in order to be able to prove that exact $H$- one forms are closed.) It is standard, going back to Driver,\cite {Driver92}, that then $d_\HH : \Dom(d_\HH )\subset L^2(C_{x_0}M\to L^2 \Gamma\cal H^*$ is closable. We will denote its closure by $\bar{d^0}$ to show it is acting on zero forms, or simply by $\bar{d}$, and let $\D^{2,1}$ be its domain with graph norm. There is the analogous result for functions with values in a separable Hilbert space $G$. In this case the domain will be written as  $\D^{2,1}(G)$ or $\D^{2,1}(C_{x_0}M;G)$ and for almost all $\sigma \in C_{x_0}M$ the derivative $\bar{d}f_\sigma$ of $f$ at the path $\sigma$ will be in the space of Hilbert-Schmidt maps $\L_2(\HH;G)$. 
As usual for real valued functions there is the corresponding  gradient operator $\nabla:\D^{2,1}\to L^2\Gamma\cal H$. The negative of its adjoint we write as $$\div:\Dom(\div)\subset L^2\Gamma \HH\to L^2(C_{x_0}M; \R),$$ so if $V$ is an $H$-vector- field in $\Dom (\div)$ and $f\in \D^{2,1}$ then 
\begin{equation}\label {eq-div-1}
\begin{split}
&\int_{C_{x_0}M}\bar{d}f(V)d\mu =\int_{C_{x_0}M} \langle \nabla(f)(\sigma),V(\sigma)\rangle _{\HH_\sigma} d\mu(\sigma) \\
&=-\int_{C_{x_0}M} f(\sigma)\; \div(V)(\sigma) d\mu(\sigma) .
\end{split}
\end{equation}
This divergence operator is closed and the standard Riesz correspondence $\phi \mapsto \phi^\Sharp$ with inverse $V\mapsto V^\Sharp$ between $H$-one-forms and $H$-vector fields maps the domain of the adjoint $d^*$ of $\bar{d}$ to that of the divergence with
$d^*\phi=-\div(\phi^{\#})$. 

For $1\leqslant p<\infty$ there are the spaces $\D^{p,1}$ defined in the same way as for $p=2$ but using $L^p$ norms. Spaces of ``weakly differentiable" functions $\mathbb{W}^{p,1}(C_{x_0}M;G)$,$1\le p<\infty$, were also given in  \cite {Elworthy-Li-Ito-map}, loosely following  \cite{Eberle-book}. Here we shall also denote those weak derivatives by $\bar{d}$. Whether $\mathbb{W}^{p,1}=\D^{p,1}$ , as occurs on $C_0\R^m$, is an open question. 
We note the following from \cite {Elworthy-Li-Ito-map}, \cf \cite {Elworthy-Li-Hodge-1}. Parts (a) and (b) are essentially equivalent and (a) is a vital step in the proof of the closability of the exterior derivative used below.
\begin {theorem}\label{th-Ito}
\begin {itemize}
\item [(a)] The map $\overline{T\I(-)}_.$ from $L^2(C_0\R^m;H)$ to vector fields on $C_{x_0}M$ given by \begin{equation}
 \overline{T\I(V)}_\sigma=\mathbb{E}\big\{\omega \mapsto T\I_\omega(V(\omega)) |\I(\omega)=\sigma\big \}
\end{equation}
gives a continuous linear map 
 $\overline{T\I(-)}_.:L^2(C_0\R^m;H)\to L^2\Gamma\cal H.$
 \item [(b)] The pull back operation $\phi\mapsto \I^*(\phi)$ defined from one- forms on $C_{x_0}M$ to $H$-one forms on $C_0\R^m$ by $\left(\I^*{\phi}\right)_{\omega}=\phi_{\I(\omega)}\circ T_\omega\I$  extends to give a continuous linear map $\I^*:L^2\Gamma \HH^* \to L^2(C_0\R^m;H^*)$.
\item[(c)] If $f\in \D^{p,1}(C_{x_0}M;G)$ then the composition  $f\circ\I$ is in $\D^{p,1}(C_0\R^m;G)$ and  then $\bar{d}(f\circ \I)=\I^*(\bar{d}f)$. 
\item[(d)] A measurable function $f:C_{x_0}M\to G$ has $f\in \mathbb{W}^{p,1}(C_{x_0}M;G)$ iff the composition  $f\circ\I$ is in 
$\D^{p,1}(C_0\R^m;G)$ and  then the weak derivative $\bar d f$ satisfies 
$\bar{d}(f\circ \I)=\I^*(\bar{d}f)$.
\end{itemize}
\end{theorem}

{\bf B. Exterior differentiation of $H$-one-forms.}
For any $C^1$ one form $\phi$ on $C_{x_0}M$ there is the usual exterior derivative $d\phi$ given by formula (\ref{Palais}). This can be restricted to give an $H$-$2$-form, $d^1_\HH $ say. As for functions we choose an initial domain $\Dom(d^1_\HH)$ to give an operator: 
$$ d^1_\HH:\Dom(d^1_\HH)\subset L^2\Gamma (\HH^1)^* \to L^2\Gamma (\HH^2)^*.$$ The domain must consist of $C^2$ one-forms $\phi$ on $C_{x_0}M$ which satisfy
\begin{enumerate}
\item [(i)] as an $H$-one-form $\phi\in L^\infty \Gamma \HH^*$. 
\item[(ii)] the exterior derivative $d\phi$ when restricted to $\HH^2$
is essentially bounded, \ie $d^1_\HH\phi \in L^\infty \Gamma {\HH^2}^*$.
\item[(iii)]  ( module structure)  if $f \in\Dom(d_\HH )$ and $\phi\in \Dom(d^1_\HH )$ then $f\phi\in \Dom(d^1_\HH )$.
\item [(iv)]  The domain of $d_\HH $ is mapped into the domain of $d^1_\HH $ by $d_\HH$.

\end{enumerate}

All these hold if we use smooth cylindrical functions and forms as initial domains, or $C^2$ functions and  $C^1$ forms which are bounded together with
their exterior derivative using the natural Finsler metric on $C_{x_0}M$. In fact it is shown in \cite {Elworthy-Li-Ito-map} that $\D^{2,1}$ is independent of the choice of $\Dom(d_\HH )$ under these restrictions, so we may as well assume that the latter is the space of smooth cylindrical functions. 

Under these assumptions we have
\begin {theorem} \cite {Elworthy-Li-Hodge-1}  The exterior derivative considered as an operator $$  d^1_\HH :\Dom(d^1_\HH )\subset L^2\Gamma (\HH ^1)^* \to L^2\Gamma (\HH ^2)^*$$ is closable. 
\end {theorem}
Since the proof was given in full in  \cite {Elworthy-Li-Hodge-1} and the analogous proof for two forms is in Part II it will be omitted here.

Let $\bar{d^1}$ denote the closure of  $d^1_\HH $. 
\begin{theorem}\cite {Elworthy-Li-Hodge-1}. The derivative $\bar{d^0}f$ of any function 
$f\in \D^{2,1}$ lies in the domain of $\bar{d^1}$ and $$ \bar{d^1}\bar{d^0}f=0.$$
\end{theorem}

The derivation property $\bar{d^1}(f\phi)=f\bar{d^1}\phi+\bar{d^0}f\wedge\phi$ is given meaning and proved in Theorem \ref{op-theorem-2} below.
\medskip

{\bf C. The first $L^2$ de Rham cohomology group and a Hodge decomposition for $H$-one-forms.} From the results above we can define the first $L^2$ cohomology group  of $C_{x_0}M$ to be the quotient of the kernel of $\bar{d^1}$ by the image of $\bar{d^0}$.
An important result here is due to Fang :

\begin {theorem} [Fang \cite{Fang94}]The image of  $\bar{d^0}$ is a closed subspace of $L^2\Gamma\HH^* $.
\end {theorem}

 It is almost a formality now to define the self adjoint Hodge-Kodaira operator $\triangle$ or
 $\triangle^1$ by  $$ \triangle^1=\bar {d^1}^* \bar{d^1}+ \bar{d^0}\bar {d^0}^*.$$
 and to obtain the Hodge decomposition. For the details we refer to \cite {Elworthy-Li-Hodge-1} or Part II.
 \begin{theorem} \cite {Elworthy-Li-Hodge-1}. There is the orthogonal decomposition
  $$ L^2\Gamma\HH = \Image (\bar{d^0}) + \overline{\Image (\bar{d^1}^*)}+ \ker\triangle^1$$ where $\overline{\Image (\bar{d^1}^*)}$ denotes the closure of the image of the adjoint of $\bar{d^1}$. 
\end{theorem}


\section{Tensor products as operators: algebraic operations on H-one forms}
\label{se-tensor}

 To show that the exterior product of $H$-one-forms
can be defined as an $H$-two-form (by a pointwise
construction) and to obtain a better understanding of the spaces
$\HH _\sigma ^2$ we will give an interpretation of $H$-two-vectors in terms of
linear maps from $\HH _\sigma ^1$ to itself. We will also give
an example on flat linear Wiener space to show how a
theory of tangent processes would lead to analogues of the elements in $%
\HH _\sigma ^2$.

\bigskip

{\bf A.}
 First we establish our notation and review the well known results identifying various completions of the algebraic tensor product $H\otimes _0H$, with spaces of linear maps, and the dualities
between the spaces. For example see Ruston \cite{Ruston}, though our conventions are slightly different. Here $H$ will be a separable real Hilbert space.
 Identify  $H\otimes _0H$ with finite rank operators on $H$ by

$$H\otimes _0H\longrightarrow \L  (H;H) $$
given by
\begin{equation}
(u\otimes v)(h)=\langle v,h\rangle u.  \label{op-1}
\end{equation}
This extends to an identification of, the projective tensor product (the
``smallest'') $H\otimes _\pi H$ with the space ${\cal L}_1(H;H)$ of
trace class operators, of our usual $H\otimes H$ with the Hilbert-Schmidt
operators ${\cal L}_2(H;H)$, and of the inductive, the `largest
reasonable', completion $H\otimes _\varepsilon H$ with the space of compact
operators ${\cal L}_c(H;H)$ in $\L  (H;H):$
\[
\begin{tabular}{lllllll}
$H\otimes _\pi H$ & \vector(1,0){50} & $H\otimes H$ & \vector(1,0){50}
&  $H\otimes _\varepsilon H$ &  &  \\
\vector(0,-1){50}  \lower0.8cm \hbox{$\simeq$}  &
  & \vector(0,-1){50} \lower0.8cm \hbox{$\simeq $}  &
& \vector(0,-1){50} \lower0.8cm \hbox{$\simeq $} &  &  \\
${\cal L}_1{\cal (}H;H)$ & $\vector(1,0){50} $ & ${\cal L}_2{\cal (}H;H)$
& \vector(1,0){50} & ${\cal L}_c{\cal (}H;H)$ &
 $\hookrightarrow $
& $\L  (H;H)$.  \end{tabular}
\]
The vertical arrows above are isometries, the inner product on
 ${\cal L}_2(H;H)$ being given by
\begin{equation}\begin{split}
\left\langle S,T\right\rangle _{{\cal L}_2} :=&
 \hbox{trace }T^{*}S  \label{op-2} \\
=&\sum_{i=1}^\infty \left\langle Se_i,Te_i\right\rangle _H
\end{split}\end{equation}
for $\left\{ e_i\right\} _{i=1}^\infty $ an orthonormal base of $H$.
So $\trace(u\otimes v)=\langle u, v\rangle$ and
$$\|u\otimes v\|_{{\cal L}_2}=\|u\|\|v\|=\|u\otimes v\|_{H\otimes H}.$$

These conventions lead to the following isomorphism with the space of 
bounded bilinear maps
$$\L  (H;H)  \longrightarrow \L  (H,H;\R ) $$
$$T \mapsto \tilde{T} $$
being given by
\begin{equation}
\tilde{T}\left( h_1,h_2\right) = \left\langle h_1,Th_2\right\rangle
\label{op-3}
\end{equation}
with resulting  isomorphism, as $\L  (H,H;\R )\simeq (H\otimes _\pi H)^*$,  
$$\L  \left( H;H\right) \stackrel{D_1}{\longrightarrow }
\left( {\cal L}_1 (H;H)\right) ^{*} $$
expressed by
\begin{equation}
D_1\left( T\right) \left( S\right) =\hbox{trace }T^{*}S.  \label{op-4}
\end{equation}

\bigskip

This construction shows that $D_1$ restricts to an isomorphism
$$\L  _{skew}(H; H) \stackrel{D_1}{\longrightarrow }
\left(\wedge_\pi^2 H\right) ^{*}  $$
where $\L  _{skew}(H;H) $ refers to the skew adjoint
elements of $\L  (H;H)$. We shall see later that  our operator $Q$
can be considered as a map from
 $\wedge^2\HH ^1$ to $\L   _{skew}(\HH;\HH)$.

\bigskip


{\bf B.}
 We will need the `double duality' map $\breve\theta=D_1^{*}\circ i$  with $i$ the canonical inclusion ${\cal L}_1(H;H)\to{\cal L}_1(H;H)^{**}$:
\begin{eqnarray*}
&& {\cal L}_1(H;H) \stackrel{_{\breve{\theta}}}
{\longrightarrow }\L  (H;H) ^{*}\\
&& \breve{\theta}(T)(S) :=\hbox{trace }S^{*}T,
\end{eqnarray*}
$T\in {\cal L}_1(H;H), \hskip 6pt S\in \L  (H;H)$.
Through the isomorphism $ {\cal L}_1(H;H) \simeq  H\otimes _\pi H $, 
it corresponds to the continuous bilinear map
$$\theta :H\times H\longrightarrow \L  \left( H;H\right) ^{*} $$
given by
$$\theta \left( h^1,h^2\right) =\breve{\theta}\left( h^1\otimes h^2\right) $$
so that
\begin{equation}
\theta \left( h^1,h^2\right) (S)=\left\langle h^1,Sh^2\right\rangle.
\label{op-5}
\end{equation}
\bigskip

{\bf C.} Let $H=L_0^{2,1}\R^m$. If $V$ belongs to the inductive
tensor product \\ ${H\otimes _\varepsilon H
    \hookrightarrow \otimes _\varepsilon ^2C_0\R^m}$
 we see, by taking $V$  primitive, that the corresponding element
  $S^V$, say, in $ \L  \left( H;H\right) $ is given by
\begin{equation}
 S^V(h)_s  \equiv V(h)_s
= \int_0^T \left(\frac  \partial {\partial
t}V_{s,t}  \right) \left(\dot{h}_t\right) \;dt, \label{op-6} \end{equation}
identifying $\frac \partial {\partial t}V_{s,t}\in \R^m\otimes \R ^m$  with the corresponding element of $\L \left(\R ^m;\R ^m\right)$. For more general kernels
$V\in  \otimes _\varepsilon ^2C_0\R ^m$ this can be used to define a linear operator $S^V$ and we
 let  ${\cal K}\R ^m$  denote the set of such $V$
for which $\frac \partial {\partial t}V_{s,t}$ exists for almost
all $t$ for each $s\in [0,T]$ and (\ref{op-6}) determines an element $S^V$
of $\L  (H;H)$.
\bigskip

As our main example of an element of ${\cal K}\R ^m$  let
$$j:[0,T]\longrightarrow \R ^m\otimes \R ^m  $$
be absolutely continuous with essentially bounded derivative
 and $j(0)=0$.
Set $V_{s,t}=j\left( s\wedge t\right) $ for $s\wedge t:=\min (s,t)$. Then
$V$ belongs to
$ {\cal K}\R ^m$ \begin{eqnarray*}
S^V\left( h\right)_s
 &=&\int_0^T {\partial \over \partial t}j(s\wedge t)(\dot h_t)dt
=\int_0^sj^{\prime }(r)(\dot{h}_r) dr
\end{eqnarray*}
and there is a conjugacy

\[\begin{tabular}{clc}
$L^2\left( [0,T];\R ^m\right) $ & $\stackrel{M^{j^\prime }}
{\vector(1,0){40}}$ & $L^2\left([0,T]; \R ^m\right)$ \\
\vector(0,1){25}  \raise0.2cm\hbox{$\frac d{d\cdot }$} &  &\vector(0,1){25}
\raise0.2cm \hbox{$\frac d{d\cdot } $} \\
$L_0^{2,1}\R ^m$ & $\stackrel{S^V}{\vector(1,0){40}}$ & $L_0^{2,1}%
\R ^m$
\end{tabular} \]\\
to the multiplication (i.e. zero order) operator $M^{j^\prime} $ given by
 $${M^{j'}(f)=j'(t)\,f(t)}$$
 for $ j^{\prime}(t)$ considered to be in $  \L  (\R ^m;\R ^m)$.
 In particular we see that in general such $V$ do not correspond to compact
  operators, let alone to elements of $H\otimes H$. Also for
${\theta: H\times H\rightarrow
\L  (H; H)^{*}}$
defined in $\S C$ we see from (\ref{op-5}) that
\begin{equation}
\theta \left( h^1,h^2\right) \left( S^V\right)
=\int_0^T\left\langle \dot{h}_s^1,\, j'(s) (\dot{h}_s^2) \right
 \rangle_{\R ^m}ds.  \label{op-7}
\end{equation}

\begin{theorem}
 For $V$ in $ {\HH _\sigma ^1\wedge {\cal
H}_\sigma ^1}$ let $ {Q(V)}\in
\wedge^2_\epsilon T_\sigma C_{x_0}$ be defined by (\ref{Q}). Then
considered as a kernel it determines  an element $S^{Q(V)}$ of
$ {\L  \left( \HH _\sigma ^1;
   \HH _\sigma ^1\right) }$
 which is conjugate to a multiplication operator
 ${M}$ on ${L^2T_\sigma C_{x_0}M}$:
\[
\begin{tabular}{clc}
$L^2T_\sigma C_{x_0}M$ & $\stackrel{M}{\vector(1,0){50}}$ &
$L^2T_\sigma C_{x_0}M$ \\
\vector(0,1){30}  \raise0.2cm\hbox{$\frac \D{\partial t}$} &
 &\vector(0,1){30} \raise0.2cm\hbox{$\frac \D{\partial t}$} \\
$\HH _\sigma ^1$ & $\stackrel{S^{Q(V)}}{\vector(1,0){50}}$
& $\HH _\sigma ^1$
\end{tabular}. \]
Here $M(u)_t=W_t {j'_V}(t) (W^{-1}_t u_t)$ for $j_V$ given by equation (\ref{eq-jay-def}) (and so $j'_V$ by
 (\ref{op-10}) below), $u\in L^2TC_{x_0}M$.
 \end{theorem}

{\it Proof.}  Set  ${\tilde V}_{s,t}= \left(W_s^{-1}\otimes W_t^{-1}\right)V_{s,t}$. 
Let
 ${\tilde{Q}:\wedge ^2L_0^{2,1}T_{x_0}M
\longrightarrow \wedge^2C_0T_{x_0}M}$
be given by
\begin{equation} \label{8.5}
\tilde{Q}(U)_{s,t}=\left( W_s^{-1}\otimes W_t^{-1}\right) Q(\wedge
^2( U)_{s,t} \end{equation}
Then  from equation (\ref{eq-jay})
$$\tilde Q({\tilde V})_{s,t}=j_V(s\wedge t).$$
As earlier $S^{\tilde{Q}({\tilde V})}$ is conjugate, by
 ${\frac d{dt}}$, to ${M^{j'_V} }$
acting on ${L^2\left( [0,T];T_{x_0}M\right)}$. 

For $h\in \HH _\sigma ^1$ we have
$S^{Q(V)}(h))_t=W_t\left(S^{ \tilde{Q}(\tilde V)}(W_{\cdot }^{-1}h_{\cdot})\right)_t$
so
\begin{eqnarray*}
\frac{{\Bbb D}}{dt}S^{Q(V)}(h)_t &=&W_t\frac d{dt}\left( S^{\tilde{Q}(\tilde V)}(W_{\cdot
}^{-1}h_{\cdot })\right)
=W_t\left( M^{j'_V} \left( \frac d{dt}W_{\cdot }^{-1}h_{\cdot }\right)
\right) _t \\
&=&W_t\left( M^{j'_V} \left( W_{\cdot }^{-1}\frac{{\Bbb D}}{dt}h_{\cdot
}\right) \right) _t
= W_t{j'_V(t)}W^{-1}_t{\D\over dt}h\\
\end{eqnarray*}
 proving the conjugacy.
\hfill\rule{2mm}{2mm}
\medskip

Thus $Q(V)_\sigma$ corresponds to an element of $\L_{\rm skew}(\HH_\sigma;\HH_\sigma)$, and so of
 $({\HH_\sigma}^*\otimes_{\pi} {\HH_\sigma}^*)^*$,
but is not compact and in particular does not belong to $\wedge^2 \HH ^1_{\sigma}$.
This yields
\begin{corollary}
\label{H2-inclusion}
There is a natural inclusion of $\HH_\sigma^2$ in $\L_{skew}(\HH_\sigma; \HH_\sigma)$ given by $V\mapsto S^V$.
\end{corollary}
Note that by the definition (\ref{eq-jay-def}) and formula (\ref{relate-R2-R})
\begin{eqnarray}
\label{op-10} 
&&j_V^{\prime }(t)\\
&=& \nonumber
\left(W_t^{-1}\otimes W_t^{-1}\right)
 (\frac{{\Bbb D}^{(2)}}{dt}+{\cal R}_{\sigma _t})
W_t^{(2)}\int_0^t\left( W_r^{(2)}\right) ^{-1}{\cal R}_{\sigma _r}
\left( \wedge ^2(W_r)V_{r,r}\right) dr 
\\ \nonumber
 &=&\left(W_t^{-1}\otimes W_t^{-1}\right)
({\cal R}_{\sigma_t}\left( \wedge ^2(W_t)V_{t,t}\right))\\
&& \nonumber
 \; +\left( W_t^{-1}\otimes W_t^{-1}\right) ({\cal R}_{\sigma
_t}W_t^{(2)}\int_0^t\left( W_r^{(2)}\right) ^{-1}{\cal R}_{\sigma
_r}\left( \wedge ^2(W_r)V_{r,r}\right) dr) .
\end{eqnarray}

\bigskip
\begin {remark} The inclusion can also be seen geometrically from the fact that if $U\in \HH^2_\sigma$ then $U-\pathR(U)\in \wedge^2\HH_\sigma\subset \L_{skew}(\HH_\sigma; \HH_\sigma) $ where $\pathR$ is the curvature operator of the damped Markovian connection which takes values in $\L_{skew}(\HH_\sigma; \HH_\sigma)$; see \S\ref{se-torsion}D below.
\end{remark}

{\bf D. Interior and exterior products.}
 For any separable Hilbert space $H$ define the interior product
 by an element  $h$ of $H$ by $$\iota_h: H\otimes_0 H\rightarrow H,
\hskip 24pt h \in H$$
$$\iota_h(h^1\otimes h^2 ):=\langle h^1, h\rangle h^2 \,\,
=S^*(h),$$
where $S\in \L  (H;H)$ corresponds to $h^1 \otimes h^2$. Thus
$\iota_h$ extends to a continuous linear map over all the completed
tensor products we use and even can be defined consistently as
$$\iota_h: \L  (H;H)\rightarrow H, \hskip 28pt \hbox{by}$$
$$\iota_h(S):=S^*(h).$$

\bigskip


{\bf E.} The first part of the following lemma is standard, but the conventions
are important, see Appendix A.

\begin{lemma}\label{op-G}
\begin{enumerate}
\item[(i)] The maps $\iota_h: H\otimes H \longrightarrow H$ and
$h\otimes: \, H\longrightarrow H\otimes H$ are mutually adjoint as are the
maps $\iota_h :\wedge^2 H \longrightarrow H$ and
 $h\wedge:\, H \longrightarrow \wedge^2 H$.

\item[(ii)] The adjoint of $h\otimes: \, H\longrightarrow H\otimes_\pi H$
is  $\iota_h: \L   (H; H) \longrightarrow H$, identifying
$(H\otimes_\pi H)^*$ with $ \L   (H; H)$ by $D_1$ as in (\ref{op-3}).
Similarly the adjoint of $h\wedge: H\to \wedge^2 _\pi H$ is the restriction
of $\iota_h$ to the skew symmetric elements $\L  _{skew}(H;H)$ of
$ \L  (H;H)$, using the restrictions of $D_1$, (see \S B above).
\end{enumerate}
\end{lemma}

{\it Proof of} (ii).  If $S\in \L  (H;H)$ and $h_1 \in H$ then
\begin{eqnarray*}
\langle\iota_h(S),h_1\rangle
&=&\langle S^*(h), h_1\rangle
=\hbox{trace }[S^*\circ (h\otimes h_1)]\\
&=& D_1 (S) (h\otimes h_1)=D_1(S)(h\otimes \cdot)(h_1)
\end{eqnarray*}
while if $S$ is skew symmetric
$$D_1(S)(h\otimes h_1)=
\langle h, Sh_1\rangle ={1\over 2} \left\{
\langle h, Sh_1\rangle -\langle Sh, h_1\rangle\right\}
= D_1 (S)(h\wedge h_1).$$
\hfill\rule{2mm}{2mm}

\bigskip
{\bf F. } Now take $H=L_0^{2,1} T_{x_0}M$ and consider
$\tilde Q:\wedge^2 H\to \wedge^2 C_0 T_{x_0}M $ given as
in (\ref{8.5}). The inclusion
${H\hookrightarrow C_0T_{x_0}M}$ has an injective adjoint 
 with dense range ${(C_0T_{x_0}M)^* \to H}$.  Let
 ${\phi^{\#}}$ denote the image of
 ${\phi \in (C_0T_{x_0}M)^*}$ under this map. There is the
 interior product
$$\iota_\phi:\wedge^2 C_0T_{x_0}M\longrightarrow C_0T_{x_0}M$$
given by
$$\iota_\phi(u^1\wedge u^2)={1\over 2}
\left(\phi(u^1)u^2-\phi(u^2)u^1\right).$$

\begin{lemma}\label{op-H}
 For $\underline h\in \wedge^2 H$ consider
 ${ S^{\tilde Q_\sigma(\underline h)}\in \L  (H;H)}$.
Then for $\phi\in (C_0T_{x_0}M)^*$
we have
$$\iota_\phi(\tilde Q_\sigma (\underline h))
=\iota_{\phi^{\#}}S^{\tilde Q_\sigma(\underline h)}
=-S^{\tilde Q_\sigma(\underline h)} (\phi^{\#}).$$
\end{lemma}
\noindent
{\it Proof.}
Write $\phi$ in terms of a $T_{x_0}M$-valued countably additive measure,
${m^\phi}$, of finite variation on $[0,T]$ so
$$\phi(w)=\int_0^T \langle w_s, dm^\phi (s)\rangle, \hskip 20pt
w\in C_0 T_{x_0}M.$$
Then, if ${\underline u=u^1\wedge u^2\in \wedge^2 C_0 T_{x_0}M}$,
\begin{eqnarray*}
\iota_\phi(\underline u)_t &=& {1\over 2} \int_0^T
 \langle u^1_s, dm^\phi(s)\rangle u^2_t
- \langle u^2_s, dm^\phi(s)\rangle u^1_t\\
&=&-\int_0^T \underline u_{t,s} \left(dm^\phi(s)\right)
\end{eqnarray*}
treating ${\underline u_{t,s} \in \wedge^2 T_{x_0}M}$ as an
element of ${\L  _{skew}(T_{x_0}M;T_{x_0}M)}$.
Thus
\begin{eqnarray}
\left(\iota_\phi \left[\tilde Q_\sigma(\underline h)\right]\right)_t
&=&-\int_0^T  j_{\underline h}(s\wedge t)\left(dm^\phi(s)\right)
\nonumber\\
&=&-\int_0^t ({d\over ds}   j_{\underline h}(s))
(\int_s^Tdm^\phi(r))ds.
\label{op-12}
\end{eqnarray}
On the other hand, if $k\in H$,
\begin{eqnarray*}
\int_0^T\langle \dot\phi^{\#}_s, \dot k_s\rangle ds
&=&\langle \phi^{\#}, k\rangle_H
=\int_0^T\left\langle k_s, dm^\phi(s)  \right\rangle ds\\
&=&\int_0^T \left\langle \dot k_s, \int_s^Tdm^\phi(r)\right\rangle ds.\\
\end{eqnarray*}
Thus ${\phi_t^{\#}=\int_0^t\left(\int_s^Tdm^\phi (r)\right)ds}$,
(a well known result in Wiener space theory). This, using (\ref{op-12}) and then \S C above, gives
\begin{eqnarray*}
\left(\iota_\phi \left[\tilde Q_\sigma(\underline h)\right]\right)_t
&=&-\int_0^t  {d\over ds}   j_{\underline h}(s)\, (\dot\phi^{\#}_s)ds
=-S^{\tilde Q_\sigma(\underline h)}(\phi^{\#}),\\
&=&\iota_{\phi^{\#}} S^{\tilde Q_\sigma(\underline h)}
\end{eqnarray*}
by definition (see \S E).
\hfill\rule{2mm}{2mm}

\bigskip

\noindent {\bf Remark.} The same calculation shows that the analogous result holds with
 general elements of ${\cal K}T_{x_0}M$, see \S C, replacing
 $\tilde Q_\sigma (h)$.
\bigskip


{\bf G. } Set
$$\tilde\HH ^2_\sigma =(1+\tilde Q_\sigma)[\wedge^2 H]
\subset \wedge^2 C_0T_{x_0}M.$$
From \S D above we can consider elements of
${\tilde\HH ^2_\sigma}$ as skew-symmetric bounded linear operators on $H$.
This can be exploited to extend the definition of exterior products:
\begin{lemma}\label{op-I}
The mapping
$$(C_0T_{x_0}M)^* \times  (C_0T_{x_0}M)^*\to (\tilde \HH ^2_\sigma)^*$$
given by
$$(\phi^1, \phi^2)\to \phi^1\wedge \phi^2|_{\tilde \HH ^2_\sigma}$$
extends to a continuous, antisymmetric, bilinear map
$$H\times H \stackrel{\wedge}{\to} (\tilde \HH ^2_\sigma)^*$$
inducing a bounded linear map
 $\tilde\theta_\sigma:
\wedge^2_\pi H\to (\tilde \HH ^2_\sigma)^*$ which
 agrees with the map $\breve\theta$ of \S C:\\
\begin{picture}(200,100)(0,0)
\put(10,10){$\wedge^2_\pi H$}
\put(35,15) {\vector(1,0){120}}
\put(95,20){$\tilde \theta_\sigma$ }
\put(170,10){$(\tilde \HH ^2_\sigma)^*$}
\put(170,80){$\L  (H;H)^*$}
\put(185,75){\vector(0,-1){55}}
\put (35,22){\vector(2,1){120}}
\put (80,45){$\breve\theta$}
\end{picture}

\noindent
using the inclusion of ${\tilde\HH ^2_\sigma}$ into $\L  (H;H)$.
\end{lemma}

\noindent{\it Proof.}
For ${S\equiv S^{\tilde Q_\sigma(\underline h)}\in \L  _{skew}(H;H)}$
 corresponding to
$\tilde Q_\sigma (\underline h)$ as above, if
 $\phi^1, \phi^2 \in (C_0T_{x_0}M)^*$ then using
Lemma \ref{op-H}.
\begin{eqnarray}
(\phi^1\wedge\phi^2)(\tilde  Q_\sigma(h))&=&
\phi^2\left(\iota_{\phi^1} (\tilde Q_\sigma(\underline h)\right)
=-\phi^2 \left(S({\phi^1}^{\#})\right) \nonumber\\
&=&-\langle {\phi^2}^{\#}, S\left({\phi^1}^{\#}\right)
\rangle _H. \label{op-13}
\end{eqnarray}
Also
\begin{equation}
\label{op-norm-of-S}
\| S\|_{\L  (H;H)}=\sup_{0\le s\le T}
|\alpha_{\underline h}(s)|
\le \hbox{const} \cdot \sup_r |h_{rr}|
\le  \hbox{const}\cdot \|{\underline h}\|_{\wedge^2 H}
\end{equation}
for ${\alpha_{\underline h}}$ the multiplication operator
corresponding to  $S$ as in \S C, i.e.\\
 ${\alpha_{\underline h}(t)={d\over dt}  j_{\underline h}(t)}$
 given by equation (\ref{op-10}). Therefore
\begin{eqnarray*}
\vert\langle {\phi^2}^{\#},  S{\phi^1}^{\#}\rangle\vert
&\le&const \, \cdot \, \|{\underline h}\|_{\wedge^2 H}\, \cdot \,
 \| {\phi^2}^{\#}\|_H\, \cdot \,    \|{\phi^1}^{\#}\| _H.
\end{eqnarray*}

This shows we have ${\tilde \theta_\sigma \in
\L  (\wedge^2 _\pi H; (\tilde\HH ^2_\sigma)^*)}$. This agrees
with $\breve\theta$, as required,  by equality (\ref{op-5}).
\hfill\rule{2mm}{2mm}

\bigskip

{\bf H.}
We now interpret these result in terms of $\HH $-forms and $\HH $
vectors on $C_{x_0}M$.

\begin{theorem}\label{op-theorem-1}
\begin{enumerate}
\item[(i)]
For $v \in \HH ^1_\sigma$ there is an interior product (annihilation
operator)
$$\iota_v: \HH ^2_\sigma\to \HH ^1_\sigma$$
which is continuous linear, and agrees with the usual $\iota_\phi$
for $\phi\in (T_\sigma C_{x_0}M)^*$ when $v=\phi^{\#}$. The map
$(v,U)\mapsto \iota_v(U)$ is in
 $\L\left(\HH ^1_\sigma, \HH ^2_\sigma;
 \HH ^1_\sigma\right) $ and is  bounded uniformly in $\sigma$.

\item[(ii)]
The map $$(T_\sigma C_{x_0}M)^* \times (T_\sigma C_{x_0}M)^*\to
(\HH ^2_\sigma)^*$$
$$(\phi^1,\phi^2)\mapsto (\phi^1\wedge \phi^2)|_{\HH ^2_\sigma}$$
extends to give a continuous linear map
$$\lambda_\sigma: (\HH ^1_\sigma)^*  \wedge_\pi (\HH ^1_\sigma)^*
\to (\HH ^2_\sigma)^*$$ which is  bounded uniformly in $\sigma$
as an element of $ \L   \left((\HH ^1_\sigma)^* \wedge_\pi
(\HH ^1_\sigma)^* ; (\HH ^2_\sigma)^*\right)$. Moreover
\item[(iii)]
 If $v\in \HH ^1_\sigma$, $\ell \in (\HH ^1_\sigma)^* $,
and $U\in \HH ^2_\sigma$
$$\lambda_\sigma(v^{\#}\wedge \ell)(U)=\ell(\iota_vU).$$
\end{enumerate}
\end{theorem}

\noindent{\it Proof.}
(i) The existence of $\iota_v$ and its properties come from Lemma \ref{op-H} and
 the bounds on $S$ noted in equation (\ref{op-norm-of-S}).

 (ii) Lemma \ref{op-I} provides the proof of (ii) with $\lambda_\sigma$ being
conjugate by ${\wedge^2 (W_\cdot)}$ to the map
 $\tilde\theta_\sigma$ of  Lemma \ref{op-I}. We see from there that
 $\tilde\theta_\sigma$ is bounded uniformly in $\sigma$ if the inclusion
 ${\HH ^2_\sigma \to \L  (H;H)}$ is. However this is
 essentially the map ${\underline h \mapsto S^{\tilde Q_\sigma(\underline h)}}$
again.

\noindent
For (iii) approximate $v^{\#}$ and $\ell$ by elements coming from
${(T_\sigma C_{x_0})^*}$. By Lemma \ref{op-H}, if
 ${U=V+Q(V)}$
$$\iota_v(U)=\iota_v(V)-S^{Q(V)}(v^{\#})$$
so
\begin{eqnarray*}
\ell(\iota_v(U))&=&\ell(\iota_v(V))-\left\langle \ell^{\#}, S^{Q(V)}(v^{\#})
\right\rangle_{\HH ^1_\sigma} \\
&=&(v^{\#}\wedge \ell)(V)+(v^{\#}\wedge \ell)\left(Q_\sigma(V)\right),
\hskip 10pt \hbox{by  }(\ref{op-13}).       \end{eqnarray*}
 \hfill\rule{2mm}{2mm}

We shall write $\lambda_\sigma(\phi\wedge\psi)$ as $\phi\wedge_\pi\psi$ when no confusion can arise.

\begin{remark}
The map $\lambda_\sigma$ is independent of the choice of the Hilbert space inner product given to  $\HH_\sigma^1$, or $\HH^2_\sigma$. Its adjoint gives a continuous map
$$\lambda^*_\sigma: \HH_\sigma^2 \to 
\big((\HH_\sigma^1)^*\wedge_\pi (\HH_\sigma^1)^*\big)^*$$
of $\HH_\sigma^2$ into the skew symmetric linear forms on $(\HH_\sigma^1)^*$.
\end{remark}
\section{The derivation property for $\overline{d^1}$}
\label{se-derivative-p}







\bigskip

{\bf A.} We can now formulate and prove the derivation property of $\overline{d^1}$. 
\begin{theorem}\label{op-theorem-2}
Suppose
${f: C_{x_0}M\to \R}$ is in ${\Dom(\bar d^0)}$ and
${\phi \in \Dom(\overline {d^1})\cap L^\infty\Gamma(\HH ^1)^*}$
with $\overline{ d^1}\phi \in L^\infty\Gamma(\HH ^2)^*$.
Then
${f\phi \in \Dom(\overline{ d^1})}$ and
$$\overline{ d^1}(f\phi)=\bar {d^0} f\wedge_\pi \phi+f (\overline{ d^1}\phi)$$
where $\wedge_\pi $ is defined above by Theorem \ref{op-theorem-1}.
\end{theorem}

\noindent{\it Proof.}
Let $\{\phi_j\}_{j=1}^\infty$ be a sequence in
$\Dom(d^1_\HH )$ with
 $\phi_j\to \phi$ in $L^2\Gamma(\HH ^1)^*$
 and  ${d^1\phi_j\to \overline{d^1}\phi}$ in
$L^2\Gamma(\HH ^2)^*$. Assume first that
 $f\in \Dom(d_\HH )$. Then $f\phi_j\to f\phi$ in
 $L^2\Gamma(\HH ^1)^*$ by the module structure of $\Dom(d^1_\HH )$, and by standard calculus
$$d(f\phi_j)=df\wedge \phi_j+f d\phi_j.$$
therefore
$$d(f\phi_j)|_{\HH ^2_\sigma}
=\lambda_\sigma\left( df|_{\HH ^1_\sigma}\wedge \phi_j|_{\HH ^1_\sigma}\right)
+f (d\phi_j) |_{\HH ^2_\sigma}$$
in the notation of Theorem \ref{op-theorem-1}.  By the uniform bound
on ${\lambda_\sigma}$ from that theorem, and taking a subsequence if necessary to assume
${\phi_j|_{\HH ^1_\sigma}\to \phi|_{\HH ^1_\sigma}}$
for almost all $\sigma$, we see
$$\lambda_\sigma\left( df|_{\HH ^1_\sigma}\wedge \phi_j|_{\HH ^1_\sigma}\right)\to \lambda_\sigma
(\bar {d^0}f_\sigma \wedge \phi_\sigma)$$
almost surely and so  in $L^2$ by the dominated convergence theorem.
Since $f(d\phi_j) \to f\overline{d^1}\phi$ and
$f\phi_j\in \Dom(d^1_\HH )$
the result follow for
$f\in \Dom(d_\HH)$.

For general $f\in \Dom(\bar d^0)$ take
 $\{f_j\}_{j=1}^\infty$ in $\Dom(d_\HH )$
with $f_j\to f$ in $L^2$ and $\bar d f_j \to \bar d f$
in  $L^2\Gamma(\HH ^1)^*$.
From above we know that $f_j \phi \in \Dom(\overline{d^1})$ with
$$\overline{d^1}(f_j\phi)
=\bar d f_j\wedge_\pi \phi -f_j(\overline{d^1}\phi), \hskip 25pt j=1 \hbox{ to }
\infty.$$
Now $\phi$ and $\overline{d^1}\phi$ are bounded so as before
 we see  $\bar d f_j \wedge_\pi \phi\to df\wedge_\pi \phi$ and
$f_j\overline{d^1}\phi\to f\overline{d^1}\phi$, both in
 $L^1\Gamma(\HH ^2)^*$, completing the proof.
\hfill\rule{2mm}{2mm}

\section{Infinitesimal rotations as divergences}
\label{se-divergences}
We will say that a $p$-vector field $V$ on $C_{x_0}M$, (or similarly on $C_0(\R^m)$),
{\it has a divergence} if there exists
$\div V\in L^1\Gamma \wedge^{p-1}TC_{x_0}M$ such that for all
smooth, bounded, cylindrical $(p-1)$-forms $\phi$ we have
\begin{equation}
\int_{C_{x_0}M}d\phi(V)\;d\mu_{x_0}
=-\int_{C_{x_0}M} \phi(\div V)\;d\mu_{x_0}.
\end{equation}
For $p=1$ from Driver\cite{Driver92} we know that not only do sufficiently
 regular elements of $L^2\Gamma\HH ^{1}$ have divergences but so do the
{\it infinitesimal rotations} $R^\alpha\in L^2\Gamma\wedge^2 TC_{x_0}M$
given by
\begin{equation}
R_t^\alpha =\parals_t\int_0^t \parals_s^{-1} \alpha_s \, dx_s
\end{equation}
where $\alpha_s: C_{x_0}M\to \L_{skew}(T_{x_s}M;T_{x_s}M)$, $0\le s\le T$, is
in $L^2$ and progressively measurable. Indeed
$$\div R^\alpha_\cdot =0.$$
For more examples of one-vector-fields with divergences see Bell \cite{Bell}, Hu-Ustunel-Zakai\cite{Hu-Ustunel-Zakai}, and Cruzeiro-Malliavin \cite{ Cruzeiro-Malliavin} 
and for $p$-vector fields see \cite {Elworthy-Li-vector-fields}. 
As in finite dimensions if  a $p$-vector field $V$ has a divergence $\div V$, when $p>1$, then $\div V$ has a vanishing divergence. In view of the looseness of the definition and the homotopical triviality of $C_{x_0}M$ we would expect that a field with a divergence which is zero would necessarily be a divergence, and we will give some evidence for this which also sheds light on the structure of our modified de-Rham complex.

First we observe that the exterior product of suitably regular $H$-vector fields in $\Dom ({\div})$  has a divergence. For this    let $V^1, V^2\in L^2\Gamma\HH ^1$. Then we have an $L^2$ section
$V^1\wedge V^2$ of $\HH ^1\wedge \HH ^1$. If $\phi$ is a smooth (bounded) cylindrical 1-form, then as discussed in Appendix B,
$$2\, d\phi(V^1\wedge V^2)=
\iota_{V^1} d\iota_{V^2}(\phi)-\iota_{V^2} d\iota_{V^1}(\phi)
-2\phi([V^1,V^2])$$
provided $V^1$, $V^2$ are sufficiently regular. Give such regularity
\begin{eqnarray*}
&&2\int_{C_{x_0}M} d\phi(V^1\wedge V^2)\;d\mu_{x_0}\\
&=&\int_{C_{x_0}M} \iota_{V^1}(\phi)\div V^2 \;d\mu_{x_0}
-\int_{C_{x_0}M} \iota_{V^2}(\phi) \div V^1\;d\mu_{x_0}
-\int_{C_{x_0}M}\phi([V^1,V^2])\;d\mu_{x_0}.
\end{eqnarray*}
Thus $V^1\wedge V^2$ has a divergence with
\begin{equation}\label{eq-div-wedge}
2\,\div(V^1\wedge V^2)=
   -(\div V^2)V^1+(\div V^1)V^2+[V^1,V^2].
\end{equation} 
The first two terms are sections of $\HH ^1$ but as is well known,
Cruzeiro-Malliavin \cite{Cruzeiro-Malliavin96} Driver\cite{Driver-Lie-bracket}, the bracket involves
 a stochastic integral of the form I for
\begin{equation}
I_t=\parals_t\int_0^t \parals_s^{-1} {\cal R}(V^1_s\wedge V^2_s) \; dx_s,
\end{equation}
i.e. an infinitesimal rotation.
The above applies in particular to $V^i=\overline{T\I}(h^i)$ for
$h^i\in W^{2,1}(C_{x_0}M;H)$, $i=1,2$. 

Also if $\underline{h}:C_{x_0}M\to \wedge^2H$ is in $W^{2,1}$,
the 2-vector field $\overline{\wedge^2 T\I}(\underline{h})$
has a divergence with $ \div \overline{\wedge^2 T\I}(\underline{h})
=\overline{T\I(\div(\underline{h}\circ\I)}$. Indeed
 for $\phi$ a smooth cylindrical one-form
\begin{eqnarray*}
&&\int_{C_{x_0}M} d\phi \big(\overline{\wedge^2 T\I}(\underline h) \big)\;d\mu_{x_0}
= \int_{C_0\R^m} \I^* (d\phi)(\underline h\circ \I) \; dP\\
&=&\int_{C_0\R^m} d( \I^*\phi)(\underline h\circ\I) \; dP
=-\int_{C_0\R^m} \I^*\phi (\div \underline h\circ\I) \; dP\\
&=&-\int_{C_0\R^m} \phi (\overline{T\I(\div (\underline h\circ\I))}) \; dP.\\
\end{eqnarray*}
Here we use the fact that since $\underline{h}\in W^{2,1}$, we have $\underline{h}\circ\I\in
\D^{2,1}\subset \Dom(\div)$. Consequently, 
\begin{equation}
\label{div-2}
\div\left(\overline{\wedge^2 T\I}(\underline h) \right)
=\overline{T\I(\div \underline h)}.
\end{equation}
(For another version of this result see \S8E.)
On the other hand
$$\overline{\wedge^2 T\I}(\underline h)
=\wedge^2 \overline{T\I}(\underline h)
+Q(\wedge^2 \overline{T\I}(\underline h)).$$
Thus:
\begin {proposition}
\label{pr-divQ} 
For $\underline {h}=h^1\wedge h^2$ with  
$h^i\in W^{2,1}(C_{x_0}M;H)$, $i=1,2$, the two-vector field
$Q(\wedge^2 \overline{T\I}(\underline h))$ has a divergence with $$\div Q(\wedge^2 \overline{T\I}(\underline h))=\overline{T\I(\div \underline h)}-\div (\wedge^2\overline{T\I}(\underline h)).$$
 \end{proposition}
Since $\overline {T\I}(\div \underline h)\in \Gamma\HH ^{1} $
we see that $\div Q(\wedge^2 \overline{T\I}(\underline h))$ must cancel out
the infinitesimal rotation term $I$  in
$\div (\wedge^2\overline{T\I}(\underline h))$.
A geometrical interpretation of this is given below, see  \S\ref{se-torsion}.
 The following result concerning
the flat Wiener space case shows how this can happen. It should be considered together with formula (\ref{eq-jay}) for $Q$ and the discussion in \S C of \S\ref{se-tensor}.

\begin{proposition}
\label{pr-rot}
Every two-vector field $V: C_0(\R^m)\to \wedge^2C_0(\R^m)$
given by $V_{s,t}=j(s\wedge t)$ for $j(t)=\int_0^t\alpha_r \; dr$,
where $\alpha_\cdot: [0,T]\times C_0(\R^m)\to\L  _{skew}(\R^m; \R^m)$
is progressively measurable with
 $\int_{C_0(\R^m)} \int_0^T|\alpha_s|^2 \;ds<\infty$, has a divergence.
It is given by
 $$\div V=\int_0^\cdot \alpha_s \; dB_s,$$
i.e. $\div V=R^\alpha$.
\end{proposition}
\noindent
{\it Proof.}
Let $f:C_0(\R^m)\to \R$ be bounded and $C^\infty$ and let
 $\ell\in C_0(\R^m)^*$.
Define the 1-form $\phi$ on $C_0(\R^m)$ by
$$\phi_\omega(v)=f(\omega)\ell(v).$$
Bounded cylindrical 1-forms can be written as sums of such forms. Then
$d\phi=df\wedge \ell$.

Let $k$ be the image of $\ell$ under the inclusions
$C_0(\R^m)^*\to L_{0}^{2,1}(\R^m)$ adjoint to the
inclusion of $L_0^{2,1}$ in $C_0$.

From (\ref{op-7}) above we see
\begin{eqnarray*}
d\phi(V)&=&
\int_0^T \langle (\dot {\overbrace{\nabla_H f}})_s,
 \alpha_s \dot k_s\rangle_{\R^m}\;ds\\
&=& df (\int_0^\cdot \alpha_s \dot k_s ds).
\end{eqnarray*}
Thus
\begin{eqnarray*}
 \int_{C_0\R^m} d\phi(V) \; dP
&=&  \int_{C_0\R^m} f(\omega) \int_0^T \langle \alpha_s \dot k_s,
 dB_s\rangle_{\R^m}\; dP(\omega)\\
&=&  -\int_{C_0\R^m} f(\omega) \int_0^T \langle  \dot k_s,
\alpha_s dB_s\rangle_{\R^m}\; dP(\omega)\\
&=&  -\int_{C_0\R^m} f(\omega)  \ell\left(\int_0^\cdot \alpha_s dB_s
\right)\; dP(\omega)\\
\end{eqnarray*}
as required. (The last equality being obvious in the (most relevant)
case when $\ell(v)=\lambda(v_{t_0})$ some $\lambda \in (\R^m)^*$,
some $0\le t_0\le T$, in which case
$\dot k_s=\chi_{[0, t_0]}(s)\lambda$).
\hfill\rule{2mm}{2mm}

\section{Differential geometry of the space $\HH ^2$ of two-vectors}
\label{se-torsion} 
 In this section we will give a bundle structure to the Bismut tangent bundle $\cal H$ and interpret the quantities $Q$ and $\mathbb{R}$ which define $\HH ^2$ in terms of a natural  connection on  $\cal H$.

\subsection*{A. The $L^2$ tangent bundle and its frame bundle}
Our Banach manifold $C_{x_0}M$ has natural structural group  
 ${C_{id}\left([0,T];O(n)\right)}$ with frame bundle
 identified with the space of paths ${C_{\pi^{-1}(x_0)}([0,T];OM)}$ in the frame bundle $OM$ of $M$, starting at any frame over $x_0$. Note that
 $C_{id}\left(O(n)\right)$ has an 
orthogonal representation on $L^2\left([0,T];\R^n\right)$, acting pointwise 
$$C_{id}\left(O(n)\right) \stackrel{\rho}{\longrightarrow}
  O\left(L^2([0,T];\R^n)\right)$$ 
$$\rho(\alpha)(f)(t)=\alpha(t)(f(t)).$$ 
For $\alpha$, $\beta$ in $C_{id}\left(O(n)\right)$,
\begin{eqnarray*} 
&&\|\rho(\alpha)-\rho(\beta)\|_{{\Bbb L}(L^2([0,T];\R^n);L^2([0,T];\R^n))} 
=\sup_{\|f\|_{\L_2}\le 1}
 \sqrt{\int_0^T \vert\alpha(s)f(s)-\beta(s)f(s)\vert^2 ds}\\ 
&&\le \sup_{\|f\|_{2}\le 1}
 \sqrt{\int_0^T  \vert f(s)\vert^2 \sup_{0\le s \le T}
|\alpha(s)-\beta(s)|^2\, ds}\\ 
&\le& \sup_{0\le s \le T} \vert\alpha(s)-\beta(s)\vert
= d(\alpha, \beta). 
\end{eqnarray*} 
Thus $\rho$ is continuous into the uniform topology and we see it is even $C^\infty$ with derivative map $T_\alpha \rho$ at $\alpha$ : 
\begin{eqnarray*}
&&T_\alpha \rho:  T_\alpha C_\id O(n)  \longrightarrow TO\big(L^2([0,T];\R^n)\big)
 \subset \L(L^2([0,T];\R^n); L^2([0,T];\R^n)) \end{eqnarray*}
 given by $T_\alpha \rho(V)(f)(t)=V(t)f(t)$.

Let  $\pi: OM\to M$ be the orthonormal frame bundle of $M$ and let
$$\tilde \pi :  C_{\pi^{-1}(x_0)}(OM)\to C_{x_0}M$$ be the principal 
$C_\id O(n)$-bundle of continuous paths $u: [0, T]\to OM$ with $\pi(u_0)=x_0$.
 From above we see that the $L^2$ tangent bundle $L^2TC_{x_0}M$ has the structure of a $C^\infty$ bundle associated  to  $C_{\pi^{-1}(x_0)}(OM)$, whose elements $u$  act as frames on it by:  
\[\begin{array}{cl} 
u:&L^2([0,T];\R^n)\longrightarrow L^2 T_\sigma C_{x_0}M, \hskip 20pt 
\sigma=\tilde \pi u\\ 
&u(f)_t=u_t(f(t)). 
\end{array}\] 
This construction determines $L^2TC_{x_0}M$ as a $C^\infty$ bundle over 
$C_{x_0}M$. It tells us what its smooth sections  (in the Fr\'echet sense) are. (For example see Remark \ref{le-nablaY} below.)

\subsection*{ B. The pointwise connection}
 Let $\tilde \nabla$ denote the \emph{ pointwise connection}  on $C_{x_0}M$, as described in greater generality by Eliasson,\cite{Eliasson67}. 
It is defined  on the  bundle $L^2TC_{x_0}M\to C_{x_0}M$ 
by 
\begin{equation}\label{pointwise-connection}
(\tilde \nabla _V U)_t={D\over ds}\left. U(\exp_{\sigma_\cdot}(sV_\cdot))_t
\right\vert_{s=0}
\end{equation}
where ${D\over ds }$ and $\exp$ come from the Levi-Civita connection 
on $TM$. Thus  
\begin{eqnarray*} 
(\tilde \nabla_VU)_t&=& 
X(\sigma_t){d\over ds}\left. \left(Y \left(\exp_{\sigma_t}(sV_t)\right) 
           U\left(\exp_{\sigma_\cdot} (sV_\cdot)\right)_t\right)
    \right\vert_{s=0}\\ 
&=& X(\sigma_t) d\left[\tilde Y(\cdot)U(\cdot)\right](V)_t, 
\end{eqnarray*} 
where the $L^2$-valued one-form  $\tilde Y:L^2T C_{x_0}M\to L^2([0,T];\R^m)$ is the lift of $Y$, i.e 
 $$\tilde Y_\sigma(V)(t)=Y_{\sigma(t)}(V(t)).$$
This says that the pointwise connection is the L-W connection in the sense of \cite{Elworthy-LeJan-Li-book}, for the lift $\tilde X$ 
of $X$ to $C_{x_0}M$.  
 
This connection is torsion free and is metric for the $L^2$ metric. 

\begin{remark}
 \label{le-nablaY}
  The pointwise derivative
$\tilde\nabla Y: TC_{x_0}M \times L^2TC_{x_0}M\to L^2([0,T];\R^m)$ is $C^\infty$. 
\end{remark} 
To see this let $\Upsilon$ be a locally defined $C^\infty$ frame field for $L^2 TC_{x_0}M$ giving a local trivialisation
over an open subset $U$ of $ C_{x_0}M$
$$\Upsilon:U\times L^2\left([0,T];\R^m\right)\to L^2 TC_{x_0}M.$$ 
Then $$\left[\tilde Y_{\sigma}\Upsilon(\sigma)(f)\right]_t
 =Y_{\sigma_t}\left(\Upsilon(\sigma)_tf(t)\right).$$ 
Its derivative is  
$$\left( \nabla_{v_t}\tilde Y\right) \Upsilon(\sigma)_tf(t)
+Y_{\sigma_t} \left(\tilde \nabla_v\Upsilon \left(f(t)\right)\right).$$

 \subsection*{ C. The bundle structure of $\HH$ and its damped Markovian connection}

Let $C_{x_0}^0M$ be a set of paths of full measure along each element of which the Levi-Civita parallel translation, $\paral$,  is defined and satisfies its basic composition properties. Then $\cal H_\sigma$ is defined for each $\sigma\in C_{x_0}^0M$ by formula (\ref{eq-H1}) with an isometry $\W_\sigma:L^2T_\sigma C_{x_0}M\to\HH _\sigma$, with inverse $\frac{\D}{d.}$. Thus we get an induced smooth vector bundle structure on $\HH ^1$, 
over $C_{x_0}^0M$ by 
$${{\Bbb D}\over ds}: 
\HH ^{1} \stackrel{\longleftarrow}{\longrightarrow} 
L^2TC_{x_0}M.$$ 

We can use this isomorphism to pull back the point-wise connection to get a metric connection ${\Nabla}$ 
on $\HH ^1$. This is the damped Markovian connection defined in a different way by Cruzeiro-Fang in \cite{Cruzeiro-Fang95, Cruzeiro-Fang}, Cruzeiro-Fang-Malliavin \cite{Cruzeiro-Fang-Malliavin-00}. The basis for a covariant Sobolev calculus using it is given in \cite{Elworthy-Li-Ito-map}. 
In particular we have a closed covariant derivative operator $\Nabla$ with domain, denoted by $\D^{2,1}\HH^1 $, in the space of $L^2$ sections of $\HH^1 $ mapping to the $L^2$ 
sections of $\L_2(\HH^1;\HH^1)$. In general we shall not distinguish between $C^0_{x_0}M$ and $C_{x_0}M$.

 Since the inverse map to ${\D\over d\cdot}$ 
is ${\mathcal W}$ it follows from equation (\ref{Ito-map}) that this connection is the L-W connection associated
to $\overline {T\I}$ in the sense of \cite{Elworthy-LeJan-Li-book}.
With this in mind define $$ \X : C_{x_0}M\times H\to \HH^1$$ 
\begin{equation}\label{Big-X} 
\X (\sigma)(h)= \overline{T\I}(h). 
\end{equation} 
 
As noted in \cite {Elworthy-Li-Ito-map} the adjoint of  $ \X $ is
 the H-valued $H$-one-form $\Y  $ given by 
$$\Y_{\sigma}(V)=\int_0^\cdot Y_{\sigma(r)}^*{\D\over dr}V_r\; dr.$$ 
This is also a right inverse to $\X$.
Suppose  that $u^1$ and $u^2$ are in  $\D^{2,1}\HH $. For $j=1,2$, set 
 $h^j(\sigma)=\Y_\sigma\big(u^j( \sigma)\big)$. Then, by  \cite{Elworthy-Li-Ito-map}, $h^j\in \D^{2,1}(C_{x_0}M; H)$ and:   
\begin{equation} 
\label{damped-1}
\begin{split}
{ \Nabla}_{u^1(\sigma)}u^2
&=\X (\sigma)\bar{d}\left[  \Y(u^2)\right](u^1(\sigma))\\
&=\X (\sigma)\bar{d} h^2 \left(\overline{T\I} (h^1(\sigma))\right)
=\X (\sigma) \left(\overline{\bar{d}(h^2\circ\I)}_{\sigma}(h^1(\sigma))\right).
\end{split}
\end {equation}
%

We saw in proposition \ref{pr-divQ} that for certain $v^1$ and $v^2$ the two-vector field $Q(v^1\wedge v^2)$ has a divergence. After the following lemma we can identify that divergence:
\begin {lemma}\label{le-adapted} Suppose $h:C_0\R^m\to H$ is adapted. Then $$ \overline {T\I h}=\overline {T\I}(\bar{h}).$$
\end{lemma}
\noindent{\it Proof.}
Set $v_t=T\I_t(h)$. Then, since $h$ is adapted we have as for equation (\ref{filter})
$$Dv_t=\nabla _{v_t}X\left( \tilde{\parals_t}d\beta _t\right) -\frac 
12\Ric^{\#}(v_t)dt+X(x_t)\dot{h}_t\,\,dt.$$ 
Now take conditional expectations as usual to get the result.
\hfill\rule{2mm}{2mm}

\begin{theorem} 
\label{th-divQ-2}
For any $\F_\star^{x_0}$ adapted vector fields $u^i\in L^p\Gamma \HH ^1$, i=1,2, some $p>2$,
\begin{equation} 
\label{div-torsion-2}
 \div Q(u^1\wedge u^2)={1\over 2}{\Bbb T}(u^1,u^2),
\end{equation} 
where ${\Bbb T}$ is the torsion of the damped Markovian connection $\Nabla$.
\end{theorem} 


 \noindent{\it Proof.} As above set  $h^j=\Y(u^j), j=1,2$. 
 Define the adapted $H$-vector fields $\tilde{h}^j,j=1,2$ on $C_0\R^m$ by $\tilde{h}^j=  h^j\circ \I$.  First assume that each $u^j$, and so $h^j$ 
 and $\tilde h^j$,  belong to $\D^{p,1}$.
 
 By the integration by parts formulae, as for the proof of (\ref{div-2}) for two-vector-fields in \S\ref{se-divergences}, and using the fact that
 $\dot{\tilde h}^j(\omega)_s  \perp \ker X(x_s(\omega))$ a.s.:  
\begin{eqnarray*}
\div( u^j)\circ\I&=&\E\left\{\div (\tilde{h^j})|\F^{x_0}\right\}=-\E\left\{\int_0^T\left \langle \dot{h}_s,dB_s\right\rangle|\F_{x_0}\right\}
\\&=&-\int_0^T \left\langle\tilde{\dot{h}}^j_s,X(x_s)dB_s\right\rangle=\div (\tilde{h}^j).
\end{eqnarray*}
In particular $\div (\tilde{h}^j)$ is $\F^{x_0}$-measurable.
Consequently, from Proposition \ref{pr-divQ} and formula~(\ref{eq-div-wedge}). 
\begin{eqnarray*} 
&&2\, \div Q( u^1\wedge u^2) 
= 2\, \overline {T\I \left(\div(\tilde h^1\wedge \tilde h^2)\right)} 
-2\, \div(u^1\wedge u^2 )\\ 
&&=\overline {T\I \left(-\tilde h^1 \div(\tilde h^2) +\tilde h^2 \div(\tilde h^1)+[\tilde h^1, \tilde h^2] \right)} 
 -(\div u^1)u^2+u^1div(u^2)-[u^1, u^2]\\ 
&&=\overline{T\I( [\tilde h^1, \tilde h^2])}-[u^1, u^2]. 
\end{eqnarray*} 
Also from (\ref{damped-1}):  
\begin{eqnarray*} 
[u^1, u^2](\sigma)&=& \X (\sigma)\left( \overline{\left(\bar{d} \;\tilde h^2\right) }_\sigma \left(h^1(\sigma)\right)-\overline{\left(\bar{d}\;  \tilde{h}^1\right)}_\sigma\left(h^2(\sigma)\right)\right)
-{\Bbb T}(u^1, u^2)(\sigma)\\ 
&=&\X (\sigma) \left(\overline{[\tilde{h}^1,\tilde{h}^2]}\right)_\sigma
-{\Bbb T}(u^1, u^2)(\sigma)\\
&=&\overline{T\I}_\sigma\left( \overline{[ \tilde h^1,  \tilde h^2]}_\sigma \right) 
  -{\Bbb T}( u^1, u^2)(\sigma)
\end{eqnarray*} 
giving
$$2 \div Q( u^1\wedge u^2) (\sigma)
=\overline{T\I( [\tilde h^1, \tilde h^2])}_\sigma -\overline{T\I}_\sigma\left( \overline{[ \tilde h^1,  \tilde h^2]}_\sigma \right) 
+{\Bbb T}( u^1, u^2)(\sigma).$$ 
For adapted vector fields the first two terms cancel by the previous lemma, so we have (\ref{div-torsion-2}) for adapted $\D^{p,1}$ vector fields.

If $u^1$, $u^2$ are adapted but not  in $\D^{p,1}$ we can choose, \cf Lemma \ref{approximate}, sequences of adapted processes $\{u^j_n\}_{n=1}^\infty$,$j=1,2$,   in $\D^{p,1}\HH$, converging to $u^1, u^2$ in $L^p$. Then as $n\to \infty$,
$$T(u^1_n, u^2_n)\to T(u^1, u^2)$$
in $L^1TC_{x_0}M$, by the formula $$\mathbb{T}(V^1,V^2)=\tilde{X}\left((\Nabla _{V^2}\tilde {Y})V^1-(\Nabla _{V^1}\tilde {Y})V^2\right)$$
given in the Appendix B.
On the other hand for any $C^\infty$ cylindrical 1-form $\phi$,
$$\int \phi \big(T(u^1_n, u^2_n)\big)=
-2\int d\phi(Q(u^1_n\wedge u^2_n))
\to-2 \int d\phi(Q(u^1\wedge u^2)).$$
 Thus for all adapted $L^p$ vector fields $u^i$, we have
$$\div Q(u^1\wedge  u^2)= 
{1\over 2} {\Bbb T}(u^1, u^2).$$ 
\hfill\rule{2mm}{2mm}

\begin{lemma}
\label{approximate}
If $u$ is an $\F^{x_0}_\star$-adapted H-vector field in $ L^p\Gamma \HH ^1$ for some $p>1$, there is a sequence $u_n \in \D^{p,1}\HH^1$ of $\F^{x_0}_\star$ adapted H-vector fields such that 
$u_n$ converges to $u$ in $L^p$.
\end{lemma}
\noindent{\it Proof.}
Set $\tilde{h}=\Y({d \over d \cdot }u)\circ\I\in L^p(C_0\R^m;L^2( [0, T]; \R^m))$. As finite chaos expansions are dense in $L^p$, let $\{\tilde{h}_n\}$ be  a sequence of functions with finite chaos expansion converging to $\tilde{h}$ in $L^p(C_0\R^m;L^2( [0, T]; \R^m))$. Define $v_n:C_{x_0}\to L^2( [0, T]; \R^m)$ by
$$(v_n\circ\I)_t=\E \{\tilde{ h}_n | \F_t^{x_0}\}. $$
Then $v_n$ belongs to $\D^{p,1}$, see \cite{Elworthy-Li-Ito-map}.  Set
$u_n=\X(\int_0^\cdot (v_n)_s ds)$ then $u_n$ converges in $L^p$ to $u$.
\hfill\rule{2mm}{2mm} 

\begin{remark}
\begin{enumerate}
\item [(1)] It is noted in Cruzeiro-Fang \cite{Cruzeiro-Fang} that the 
divergence of $ {\Bbb T}(v^1,v^2)$ vanishes for a class of adapted H-vector 
fields $v^1$ and $v^2$.
\item[(2)] The conclusion of the the theorem does not hold for general smooth nonadapted 
vector fields. In fact for  a smooth, cylindrical,  $f: C_{x_0}M\to\R$ we have
${\Bbb T}(f\bar v^1, \bar v^2)=f{\Bbb T}(\bar v^1, \bar v^2)$. But 
\begin{eqnarray*} 
\div Q((f\bar v^1)\wedge \bar v^2) 
=\div \left(f Q(\bar v^1\wedge \bar v^2)\right) 
= f\,\div \left (Q(\bar v^1\wedge \bar v^2)\right) 
+\iota_{\nabla f}Q(\bar v^1\wedge \bar v^2). 
\end{eqnarray*} 
\end{enumerate}
\end{remark}


Though we state the following for Brownian motion measures and the damped Markovian connections note that it applies in considerable generality, for example for any metric connection on a finite dimensional Riemannian manifold with smooth measure. In it we consider the closed covariant derivative operator 
$${\Nabla}:\D^{2,1}\subset L^2 \Gamma \HH ^1 \to L^2\Gamma{\L}_2(\HH ^1; \HH ^1)$$ with $L^2$-adjoint $\Nabla^*: L^2\Gamma{\L}( \HH ^1; \HH ^1) 
 \to L^2\Gamma {\mathcal H}^1$.

\begin{proposition} 
\label{pr:divergence-of-two-tensors} 
Let $U, V\in L^\infty \Gamma \HH ^1$. 
 Suppose $U\in \D^{2,1}$ 
and ${V\in \Dom(\div)}$. 
Then ${\sigma \mapsto U(\sigma)\otimes V(\sigma)}$ as 
an element, $U\otimes V$, of 
 ${L^2\Gamma(\HH ^1\otimes \HH ^1)}$ 
is in ${\Dom({\Nabla}^*)}$ and 
\begin{equation} 
\label{divergence-of-two-tensors} 
{\Nabla}^*(U\otimes V)(\sigma) 
=-(\div V)(\sigma) U(\sigma) -{\Nabla}_{V(\sigma)}U. 
\end{equation}
In particular this holds if $U$ and $V$ are both essentially bounded and in $\D^{2,1}$ in which case:
\begin{equation}
{\Nabla}^*(U\wedge V)=\div(U\wedge V)+\frac{1}{2}{\Bbb T}\left(U,V\right).
\end{equation} 
\end{proposition} 
\noindent{\it Proof.}
 Let $Z\in \D^{2,1}\HH ^1$. 
 By 
(\ref{op-1}) and (\ref{op-2}),
\begin{eqnarray*} 
&&\int_{C_{x_0}M}\big \langle  ({\Nabla}Z)_\sigma, 
 U\otimes V(\sigma) \big\rangle_ 
 {{\cal  H}^1_\sigma\otimes \HH ^1_\sigma} \; d\mu_{x_0}(\sigma)\\ 
&&=\int_{C_{x_0}M} \left\langle  ({\Nabla}Z)_\sigma, 
  U\otimes V(\sigma) \right\rangle_ {{\cal L}_2(\HH^1 ;\HH ^1)} \; 
  d\mu_{x_0}(\sigma)\\ 
&&=\int_{C_{x_0}M}  \sum_{i=1}^\infty \left \langle 
 ({\Nabla}_{e_i}Z)_\sigma, U(\sigma)\langle V(\sigma),e_i\rangle 
 \right\rangle_{\HH ^1_\sigma }\; d\mu_{x_0}(\sigma)\\ 
&&=\int_{C_{x_0}M} \left\langle 
 {\Nabla}_{V(\sigma)}Z, U(\sigma)\right\rangle_{\HH ^1_{\sigma}} 
  \; d\mu_{x_0}(\sigma) \\ 
&&=\int_{C_{x_0}M} d\left\langle Z, U\right\rangle_{\HH ^1} (V(\sigma)) 
  \; d\mu_{x_0}(\sigma) 
-\int_{C_{x_0}M}\left\langle Z, {\Nabla}_{V(\sigma)}  U\right\rangle 
 _{\HH ^1_\sigma} d\mu_{x_0}(\sigma), 
\end{eqnarray*} 
since ${\Nabla}$ is a metric connection. This proves the first part.

 For the second part first note from \cite{Elworthy-Li-Ito-map} that $H$-vector fields  which are are in $\D^{2,1}$ are in $\Dom (\div)$. Then  plug $U\wedge V={1\over 2}\left\{U\otimes V-V\otimes U\right\}$ into equation (\ref{divergence-of-two-tensors}) and use formula (\ref{eq-div-wedge}) to see: 
\begin{eqnarray*} 
{\Nabla}^*(U\wedge V) 
&=&{1\over 2}\left\{ -(\div V)U+(\div U)V-{\Nabla}_VU+\Nabla_UV 
   \right\}\\ 
&=&\div(U\wedge V)+\frac{1}{2}{\Bbb T}(U,V).
\end{eqnarray*}
\hfill\rule{2mm}{2mm}

By formula (\ref{div-torsion-2}) this immediately gives:
\begin{corollary} 
\label{cor-nabla-star}
For $U$, $V$ as in Proposition \ref{pr:divergence-of-two-tensors} 
\begin{equation} 
{\Nabla}^*(U\wedge V)=\div(I+Q)(U\wedge V) 
\end{equation} 
provided $U,V$ are non-anticipating. In particular for $h^1, h^2$ 
in $L_0^{2,1}(\R^m)$ non-random 
\begin{equation}\div\left(\overline{\wedge^2 T\I}(h^1\wedge h^2)\right) 
=\nabla^*\left(\overline{T\I}(h^1)\wedge 
\overline{T\I}(h^2)\right). 
\end{equation} 
\end{corollary} 
%
\hfill\rule{2mm}{2mm}
 \medskip
 
Note that for $Z$ as above, if $f:C_{x_0}M\to \R$ is smooth and cylindrical then 
\begin{eqnarray*} 
&&\int_{C_{x_0}M} 
\left\langle {\Nabla}Z, f U\wedge V\right\rangle 
  _{\HH ^1\otimes \HH ^1}\; d\mu_{x_0}\\ 
&&=\int_{C_{x_0}M} 
\left\langle {\Nabla}(fZ)-Z\otimes{\nabla f}, U\wedge V\right\rangle 
  _{\HH ^1\otimes \HH ^1}\; d\mu_{x_0}\\ 
&&=\int_{C_{x_0}M}\left\{ 
\left\langle Z, f \Nabla^*(U\wedge V)\right\rangle 
-{1\over 2} \left\langle Z, U\right\rangle df(V) 
+{1\over 2}\left\langle Z,V\right\rangle df(U) \right\}\; d\mu_{x_0}.\\ 
\end{eqnarray*} 
So 
\begin{eqnarray*}\Nabla^*[f U\wedge V]&=&f \Nabla^*(U\wedge V) 
   -{1\over 2}\left\{U df(V)-Vdf(U)\right\}\\
&=&  f \Nabla^*(U\wedge V) 
   +\iota_{\nabla f}(U\wedge V)
   \end{eqnarray*}
whereas 
$$\div(I+Q)(fU\wedge V) 
=f \div(I+Q)(U\wedge V)+\iota_{\nabla f}(U\wedge V)
+\iota_{df}Q(U\wedge V).$$ 
Thus the formula is not true, if `non-anticipating' is dropped. 
 
\subsection*{D. The curvature operator} 
 
The curvature operator $\pathR$ of the damped Markovian connection 
$\Nabla$ on $\Gamma{\mathcal H}^1$ is conjugate to the curvature 
operator $\tilde{\mathcal R}: \wedge^2 TC_{x_0}M\to {\L}_\emph{skew}(L^2TC_{x_0}M;L^2TC_{x_0}M)$ of the pointwise
connection on the $L^2$ tangent bundle via the map ${\D \over dt}$. In 
fact $$\pathR: \wedge^2T_\sigma C_{x_0}M \to { 
\L}_\emph{skew}(\HH^1_\sigma; \HH^1_\sigma)$$ is given by 
 
$$(\pathR(U)h)_t=\W_t(\tilde{\mathcal R}_{\sigma}(U(\sigma))({\D\over 
d\cdot}h_\cdot)), $$ 
that is 
\begin{equation} 
(\pathR(U)h)_t 
=W_t\int_0^t W_s^{-1} {\mathcal R}_{\sigma_s}(U_{s,s}) ({\D\over 
d\cdot}h_s)\; ds. 
\end{equation} 
We shall show that this agrees with the definition given in equation (\ref{eq:R}).
 
 Our convention that $(a\otimes b)(u)=\langle b, u\rangle 
a$  makes  clear  the correspondence between the curvature operator $\RR$  of $M$ considered as a map $\RR: 
\wedge^2 TM\to {\L}(TM; TM)$ and it considered as a map $\mathcal R: \wedge^2 
TM \to \wedge^2 TM$. Note also that for a linear map $T$ 
$$[(T\otimes \1)(a\otimes b)](u)=T((a\otimes b)(u)).$$ 
Then 
\begin{eqnarray*} 
\pathR(U)(h)_t &=&W_t\int_0^t W_r^{-1} \RR(U_{rr})({\D\over 
dr}h_r)\;dr\\ 
&=&W_t\int_0^t \left[W_r^{-1}\otimes \1)\RR(U_{rr})\right] 
({\D\over dr}h_r)\;dr\\ 
&=&\int_0^t \left[W_t(W_r)^{-1}\otimes \1)\RR(U_{rr})\right] 
({\D\over 
dr}h_r)\;dr\\ 
&=&\int_0^T \chi_{[0,t)}(r)(W_t\otimes W_r)\wedge^2(W_r^{-1}) 
\RR(U_{rr}) ({\D\over dr}h_r)\;dr\\ 
\end{eqnarray*}

\begin{proposition} 
\label{pr-path-curvature} As a linear map from $\wedge^2 T_\sigma 
C_{x_0}M$ to $\wedge^2 TC_{x_0}M$, the curvature operator of the
damped Markovian connection on $\HH^1$ is given by:
\begin{equation} 
\pathR(U)_{s,t}=(W_s\otimes W_t)\int_0^t \wedge^2(W_r)^{-1}\RR(U_{rr})dr, 
\hskip 18pt t<s. 
\end{equation} 
\end{proposition} 
 
\noindent {\it Proof.} Since $\pathR(U)$ is regular, its integral representation is
$$\pathR(U)(h_\cdot)_t=\int_0^T(\1\otimes {\D\over 
dr})\pathR(U)_{t,r}({\D h_r\over d r})dr.$$ 
 Compare this with  the integral representation above the 
proposition to see the result. \hfill\rule{2mm}{2mm} 

\subsection*{E. The domain of $\bar{d^1}^*$}

An important result for functions on $C_0\R^m$ was that the domain of the divergence acting on $H$-vector fields contains $\D^{2,1}(C_0\R^m;H)$, in particular $H$-vector fields which are in $\D^{2,1}$ are Skorohod integrable, \cite {Kree-Kree}. For $C_{x_0}M$ the analogous result was proved in \cite{Elworthy-Li-Ito-map} using the damped Markovian connection. We have not yet given a "bundle structure" or connection to $\cal {H}^2$ or its dual, but $\wedge^2L^2TC_{x_0}M$ is a smooth bundle and inherits a connection from the pointwise connection. This will be the LW connection for $\wedge^2\tilde{X}$. As discussed, in general, in \cite{Elworthy-Li-Ito-map} a section $Z$ of $\wedge^2L^2TC_{x_0}M$ is in $\D^{2,1} \wedge^2L^2TC_{x_0}M$ if $\wedge^2\tilde{Y}(Z)$ is in $\D^{2,1}\left(C_{x_0}M;\wedge^2 L^2([0,T];\R^m)\right)$.  Where defined,  the map 
$$ (\1+Q)\wedge^2 \W:\wedge^2L^2TC_{x_0}M\to \HH ^2$$
is an isometry and it would be natural to use this to give a connection on $\HH ^2$.
In this sense the following shows that the results mentioned above extend to our $H$-two-forms (or equivalently for the divergence operator on $H$-two-vectors). It is stated in terms of the weak Sobolev class $W^{2,1}$ for, possibly, greater generality.

\begin {theorem} \begin{enumerate}
\item Let $\phi\in L^2\Gamma \HH ^2$. If $$\phi\circ\left(\1+Q\right)\circ\wedge^2\W\in W^{2,1}\Gamma\wedge^2\left(L^2TC_{x_0}M\right)^*$$ 
then $\phi\in \Dom(\bar{d^1}^*)$.
\item
More generally $\phi\in \Dom(\bar{d^1}^*)$ if the conditional expectation of its pull back by the Ito map $$\E\{\I^*(\phi)|\F^{x_0}\}:C_0\R^m\to (\wedge^2 H)^*$$ is in the domain of $ \bar{d^1}^*$ on $C_0\R^m$. If so, for almost all $\sigma\in C_{x_0}M$ the H-vector field $\div \phi^{\#}$ is given by
$$ \div (\phi^{\#})= \overline{T\I(\div (\E\{\I^*(\phi)|\F^{x_0}\})^{\#})}.$$

\end{enumerate}
\end {theorem}
\noindent{\it Proof.}
Set
$$g(\sigma)=\phi\circ\left(\1+Q\right)\circ\wedge^2\W\circ \wedge^2\tilde {X}(\sigma)\circ \wedge^2 (\frac{d}{d\cdot}) $$
 
 for $\sigma\in C_{x_0}M$. Then our first condition implies that $g\in W^{2,1}\left(C_{x_0}M; \wedge^2 H\right)$. Note that $g= \phi\circ\left(\1+Q\right)\circ\wedge^2\W\circ \wedge^2\Bbb X $ and so $$g\circ\I=\E\{\I^*(\phi)|\F^{x_0}\}.$$
By \cite {Elworthy-Li-Ito-map} $g\circ\I\in \D^{2,1}$ on $C_0\R^m$. By \cite {Shigekawa-Hodge} this implies that as an $H$-two-form $g\circ\I$ is in the domain of ${d^1}^*$.
 Now let $\psi\in \Dom (d^1_\HH )$, cylindrical one-form on $C_{x_0}M$. Then we have:
 \begin {eqnarray*}
 \int_{C_{x_0}M}\langle d^1_\HH\psi,\phi\rangle_{{\HH^2}^*}&=&\int_{C_{x_0}}\langle d^1_\HH \psi\left(\overline{\wedge ^2T\I}(-)\right),\phi \left(\overline{\wedge^2 T\I}(-)\right)\rangle_{({\wedge^2 H})^*}\\ 
 &=& \int_{C_0\R^m}\langle d^1_\HH \psi\wedge ^2T\I(-),\E\{\I^*(\phi)|\F^{x_0}\}\rangle_{{\wedge^2 H}^*}\\
 &=& \int_{C_0\R^m}\langle \I^*(d^1_\HH \psi),\E\{\I^*(\phi)|\F^{x_0}\}\rangle_{{\wedge^2 H}^*}\\
 &=&\int_{C_0\R^m}\langle\bar{d^1}\I^*(\psi),\E\{\I^*(\phi)|\F^{x_0}\}\rangle_{{\wedge^2 H}^*}\\
&=&\int_{C{_0}\R^m}\langle\ I^*(\psi),(\bar{d^1})^*\E\{\I^*(\phi)|\F^{x_0}\}\rangle_{{\wedge^2 H}^*}.
 \end{eqnarray*}
 
 From this the results follow.
\hfill\rule{2mm}{2mm}

\begin{corollary} Every $C^1$ cylindrical 2-form on $C_{x_0}M$ is in the domain of $\bar{d^1}^*$.
\end {corollary}
{\it Proof} Let $M^{(k)}=M\times M...\times M$ be the Cartesian product of $k$ copies of $M$ and for $0\le t_1\le \dots \le t_k\le T$ define $\rho_{\underline{t}}: C_{x_0}M\to M^{(k)}$ by $\rho_{\underline{t}}(\sigma)=(\sigma(t_1), \dots,\sigma(t_k))$.
Suppose $\phi=\rho_{\underline{t}}^*(\theta)$ for $\theta$ a $C^1$ two -form on $M^{(k)}$. Then \begin{eqnarray*}\E\{\I^*(\phi)|\F^{x_0}\}&=&\phi_{\I(\cdot)}\circ (\1+Q_{\I(\cdot)})\circ \wedge^2\Bbb{X}( \I(\cdot))\\
&=&\theta\circ \wedge^2 X^{(k)}(\I(\rho_{\underline{t}(\cdot)}))\circ \wedge ^2 Y^{(k)}\wedge^2 T\rho_{\underline{t}}\circ (\1+Q_{\I(\cdot)})\circ \wedge^2\Bbb{X}( \I(\cdot))
\end{eqnarray*}
where $X^{(k)}(z_1,...,z_k)=\oplus_{j=1}^kX(z_j):\oplus^k\R^m\to T_zM^{(k)}$ for $z=(z_1,...,z_k)\in M^{(k)}$.
Now, from the differentiability of $\theta$ and $X$ it is clear that $\theta\circ \wedge^2 X^{(k)}(\I(\rho_{\underline{t}(\cdot)}))$ is in $\D^{p,1}$ for all $1\leqslant p<\infty$, while it follows from standard approximation arguments that so is $\wedge ^2 Y^{(k)}\wedge^2 T\rho_{\underline{t}}\circ (\1+Q_{\I(\cdot)})\circ \wedge^2\Bbb{X}( \I(\cdot))$, for example as in \cite {Aida97}. Thus we can apply the theorem as required.
 \hfill\rule{2mm}{2mm}

\section{Appendix A. Conventions} In the past we have used
different conventions on the exterior product of a differential
form, inner product of two antisymmetric tensor vectors, and the
interior product of a vector with another.
Here we were driven by the insistence that exterior product spaces are subspaces of the corresponding tensor products. To make these
differences more transparent and easier for the reader to compare
to their own conventions, we list in this section the conventions we
use. It is only necessary to state them for uncompleted tensor products.

{\bf A.}
Let $E, F$ be a real linear spaces. Any
 multilinear $\psi: E\times E\times...\times E\to F$ determines a linear map 
 $\tilde{\psi}:E\otimes_0 E\otimes_0...\otimes_0 E\to F$ with $$ \tilde{\psi}(u_1\otimes\dots\otimes u_q)=\psi(u_1,...,u_q).$$

Denote by $\wedge^q_0E$ the subspace of anti-symmetric
tensors and use the convention:
\begin{equation}
u_1\wedge\dots\wedge u_q={1\over q!}\sum_\pi (-1)^\pi
u_{\pi(1)}\otimes\dots\otimes u_{\pi(q)}
\end{equation}
where the summation is over all permutations $\pi$ of $\{1,2\dots, q\}$
and $(-1)^{\pi}$ is the sign of the permutation. This convention ensures that if $\psi$ is anti-symmetric then $$\tilde{\psi}(u_1\wedge\dots\wedge u_q)=\psi(u_1,...,u_q).$$

 An inner product $\langle-,-\rangle$  on $E$ determines one on the tensor products :
 \begin{equation}
\langle u_1\otimes\dots\otimes u_q,
v_1\otimes\dots\otimes v_q\rangle =\Pi_{i=1}^q \langle u_i,v_i\rangle,
\end{equation}
any $u_i, v_i\in E$. In turn this determines one on the exterior powers by restriction, giving:
\begin{equation}
\langle u_1\wedge\dots\wedge u_q, v_1\wedge \dots\wedge_q\rangle
={1\over q!} \hbox{det } \left(
\begin{array}{cccc}
\langle u_1, v_1\rangle &\langle u_1, v_2\rangle &\dots,
 &\langle u_1, v_q\rangle\\
\dots & \dots & \dots & \dots\\
\langle u_q, v_1\rangle &\langle u_q, v_2\rangle &\dots,
 &\langle u_q, v_q\rangle\\
\end{array}
\right).
\end{equation}
Now suppose there is a pairing $\ll-,-\gg:E'\times E\to \R$ between $E$ and a linear space $E'$. We are thinking of the cases $E=E'$ with inner product pairing or $E'$ being the dual space of $E$ with respect to some topology, with $\ll l,e\gg=l(e)$. Then if $l\in E'$, there is the standard interior product, or annihilation operator $\iota_l$,

\begin{equation}\iota_l \left(u_1\otimes \dots \otimes u_q\right)
=\ll l,u_1\gg \left(u_2\otimes \dots \otimes u_q\right)
\end{equation}
This gives
\begin{equation}
\iota_l (u_1\wedge \dots \wedge u_q)={1\over q} \sum_{j=1}^q
(-1)^{j+1}\ll l,u_j\gg
u_1\wedge \dots \wedge \widehat{u_j}\wedge \dots \wedge u_q
\end{equation}
where $\widehat{u}$ means the omission of the vector $u$.
Note that \begin {itemize}
\item[{(i)}] If $E=E'$ with inner product pairing then for each $v\in E$ the operator
$\iota_v:\wedge^q_0E\to\wedge^{q-1}_0E$ is adjoint to the map determined by 
$ u_1\wedge\dots \wedge u_{q-1}\to v\wedge u_1\wedge\dots \wedge u_{q-1}$
\item[{(ii)}] The interior product is now not a skew-derivation, \cf \cite{Kobayashi-Nomizu-I},  p65. For example if $q=2$ we have $$ \iota_l(u_1\wedge u_2)=\frac{1}{2}\bigg\{\ll l,u_1\gg u_2-\ll l,u_2\gg u_1\bigg\}$$
\end{itemize}
Keeping the duality between the interior product and the ``creation operator" $v\wedge -$, for $\psi$ as above and $v\in E$ define :$$ \iota _v {\psi}: \textsf{X}^{(q-1)}E\to \R$$ by 
$$\iota_v(\psi)(u_1,\dots,u_{q-1})=\psi(v,u_1,\dots,u_{q-1}),$$
 so that if $\psi$ is skew-symmetric we have $$\iota_v(\psi)(u_1\wedge\dots\wedge u_{q-1})=\psi(v\wedge u_1\wedge\dots\wedge u_{q-1}).$$
If $\phi_1$ and $\phi_2$ are in a dual space to $E$ then $\phi_1\wedge\phi_2$ is defined on $\wedge^2_0 E$ by 
$$ \phi_1\wedge\phi_2 (u_1\wedge u_2)=\frac{1}{2}\left[\phi_1(u_1)\phi_2(u_2)-\phi_2(u_1)\phi_1(u_2)\right].$$
This is in agreement with $\iota_v(\phi_1\wedge\phi_2):=\frac{1}{2}(\phi_1(v)\phi_2-\phi_2(v) \phi_1)$.

\bigskip
 {\bf B.} More generally if $S: E_1\to E_2$ and
$T:F_1\to F_2$ are two linear maps of Banach spaces,  there is
the induced linear map $$S\otimes T: E_1\otimes_0 F_1\to E_2\otimes_0 F_2.$$
If $E_1=F_1$ and $E_2=F_2$ set $S\wedge T={1\over 2}(S\otimes T+T\otimes S$)
 so $S\otimes S$ agrees with  $S\wedge S$ as an linear operator on
 $\wedge^2E_1$. This reduces to the previous definitions when $E_2=F_2=\R$ after identifying $\R\otimes\R$ with $\R$.

\bigskip

{\bf C.} Consider now the tangent bundle $TM$ of a smooth manifold $M$.
The exterior differentiation $d: \wedge^qTM\to \wedge^{q+1}TM$ is defined
by: 
\begin{equation}
\begin{array}{ll}
&d\phi\left(V^1\wedge\dots \wedge V^{q+1}\right)\\
&={1\over (q+1)}\sum_{i=1}^{q+1} (-1)^{i+1} L_{V^i}\left[\phi
\left(V^1\wedge\dots\wedge \widehat{V^i} \wedge\dots \wedge V^{q+1}
\right)\right]\\
&+{1\over(q+1)}{\sum_{1\le i<j\le q+1} (-1)^{i+j}  }
 \phi\left( [V^i,V^j]\wedge  V^1\wedge\dots
  \widehat{V^i}\wedge \dots \widehat{V^j}\dots
 \wedge V^{q+1}\right)
\end{array}
\end{equation}
where
 $L_{V^i}$ denotes Lie
differentiation in the direction of $v^i$. 
%
This differs from the convention used in our previous research paper,  \eg \cite{Elworthy-LeJan-Li-book}, \cite{Elworthy-Li-forms-CR} \cite{Elworthy-Li-vector-fields} where we did not add any constants
before $d$ and $d^*$. 
This lead to a
change in the divergence of $q$-vector fields by a factor of $q$
\begin{equation}
\div_{old}(V)=q\; \div_{new}(V).
\end{equation}

By our conventions if $f$ is a function on $M$,
\begin{eqnarray}
\langle \langle df\wedge\phi, \psi\rangle
&=& \langle \phi, \iota_{df}\psi \rangle.\\
d(f\phi)&=&df\wedge \phi+ f\, d\phi,\\
div(fV)&=&f(div V)+\iota_{V}(df).
\end{eqnarray}

\section{Appendix B: Brackets of vector fields, torsion, and $d\phi(v^1\wedge v^2)$}
Lie brackets of H-vector fields have been discussed in many places, e.g. \cite{Driver-Lie-bracket}, \cite{Cruzeiro-Malliavin}, \cite{Leandre-cohomology96}, for completeness, and definitiveness, we give a definition and some properties here. The torsion of the damped Markovian connection is also described, for explicit formulae see \cite {Cruzeiro-Fang}. 
We refer to \cite{Elworthy-Li-Ito-map} for the  Sobolev calculus of sections of $\HH$, related bundles, and smooth bundles such as $L^2TC_{x_0}M$. The latter will always be taken here with its pointwise connection.

\begin {proposition} The inclusion map of $\HH$ into $L^2TC_{x_0}M$ is in $\D^{p,1}$ for $1\leqslant p<\infty$ as a section of $\L_2(\HH; L^2TC_{x_0}M)$ and any H-vector field $V$ in $D^{p,1}\HH$, or $\mathbb{W}^{p,1}\HH$, is a $\D^{p,1}$, or $ \mathbb{W}^{p,1}$, section of $L^2TC_{x_0}M$.
Moreover for such $V$ the pointwise (weak) covariant derivative $\tilde{\nabla}_-V$ is an $L^p$ section of $\L(\HH;TC_{x_0}M)$.  
\end{proposition}
\noindent{\it Proof.}
For the first assertion it suffices to show that the map $$\Theta: C_{x_0}M\to \L_2(H; L^2([0,T];\R^m))$$ given by
$$ \Theta(\sigma)(h)=\tilde{Y}_\sigma \overline{T\I}_\sigma(h)$$ is in $\D^{p,1}$.
However $ \Theta(\sigma)(h)_t=Y_{\sigma(t)} W_t\int_0^tW_s^{-1}X(\sigma(s))(\dot{h}_s)ds$ and so the result holds from standard arguments, as in \cite{Aida97}. For the claim about sections we can apply the corresponding arguments to 
$\sigma\mapsto \Theta(\sigma)(U(\sigma))$ for $U\in \D^{p,1}(C_{x_0}M;H)$, or in $\mathbb{W}^{p,1}(C_{x_0}M;H)$; in the latter  case it is only necessary to consider the composition with $\I$, see Theorem \ref{th-Ito}. In particular the final assertion comes from  standard results giving the continuity in $t$ of the derivative of $(\Theta \circ\I)(U\circ\I)_t:C_0\R^m\to \R^m$ eg as \cite {Nualart-book} page 106. Alternatively the derivative can be calculated explicitly as in \cite{Aida97}.
\hfill\rule{2mm}{2mm}

\begin {definition}. If $V^1$ and $V^2$ are in $\mathbb{W}^{p,1}\HH$ define their Lie bracket by 
$$[V^1,V^2]=\tilde{\nabla}_{V^1}V^2-\tilde{\nabla}_{V^2}V^1,$$
where $\tilde \nabla$ is the pointwise connection defined by formula (\ref{pointwise-connection}).
\end {definition}
By the Proposition $[V^1,V^2]$ is then a measurable vector field, i.e. section of $TC_{x_0}M$. Since the pointwise connection restricts to a torsion free connection on 
$TC_{x_0}M$ this definition agrees with the usual one. Moreover if $f:C_{x_0}M\to \R$
is smooth and cylindrical we have 
$$\bar{d}(\bar{d}f(V^2))V^1=\tilde{\nabla}_{V^1}(\bar{d}f)V^2+\bar{d}f(\tilde{\nabla}_{V^1}V^2)$$ so that $$\bar{d}(\bar{d}f(V^2))V^1-\bar{d}(\bar{d}f(V^1))V^2=\bar{d}f([V^1,V^2])$$
as usual. The torsion $\mathbb{T}(V^1,V^2)$ is defined as a measurable vector field by 
$$\mathbb{T}(V^1,V^2)=\Nabla_{V^1}V^2-\Nabla_{V^2}V^1-[V^1,V^2].$$

To see the torsion as an ``H-tensor field" use the LW characterisation of the pointwise connection to observe first that \begin{eqnarray*}
\mathbb{T}(V^1,V^2)&=&
\Nabla_{V^1}V^2-\Nabla_{V^2}V^1-\tilde{\nabla}_{V^1}V^2+\tilde{\nabla}_{V^2}V^1
\\&=& \Nabla_{V^1}V^2-\tilde{X}\bar{d} (\tilde{Y}V^2)V^1- \Nabla_{V^2}V^1+\tilde{X} (\bar{d}(\tilde{Y}V^1)V^2.
\end{eqnarray*}
Now consider the restriction  of $\tilde{Y}$ to $\HH$ as a section of $\L_2(\HH; L^2([0,T];\R^m))$. As above it lies in $\D^{p,1}$, $1\leqslant p <\infty $ with $\Nabla\tilde {Y}$ a section of $\L_2(\HH;\L_2(\HH; L^2([0,T];\R^m))$. Then
$$\Nabla_{V^1}V^2-\tilde{X}\bar{d} (\tilde{Y}V^2)V^1=-\tilde{X}(\Nabla_{V^1}\tilde{Y})V^2$$
 and so 
$$\mathbb{T}(V^1,V^2)=\tilde{X}\left((\Nabla _{V^2}\tilde {Y})V^1-(\Nabla _{V^1}\tilde {Y})V^2\right).$$
From this we see we can consider the torsion as a section of $\L_2(\wedge^2\HH; TC_{x_0}M)$. Alternatively noting that $\tilde{Y}$ maps $\HH$ into $C_0([0,T];\R^m)$ and arguing as before we see that it gives a section of $\L_{skew}(\HH,\HH;TC_{x_0}M)$. In both cases the sections are in $L^p$ for all $1\le p<\infty$.

Finally we give the result used in \S\ref{se-divergences}.
\begin{proposition} If $\phi$ is a smooth cylindrical 1-form and $V^1$, $V^2$ are in $ \mathbb{W}^{p,1}\HH $ then, almost surely,
$$2d\phi(V^1\wedge V^2)= \iota_{V^1}\bar{d}\iota_{V^2}\phi-\iota_{V^2}\bar{d}\iota_{V^1}\phi -\phi([V^1,V^2]).$$
\end{proposition}
\noindent{\it Proof.}  Using the pointwise connection on the sections of $T^*C_{x_0}M$ :
\begin{eqnarray*}
&& \iota_{V^1}\bar{d}\iota_{V^2}\phi-\iota_{V^2}\bar{d}\iota_{V^1}\phi -\phi([V^1,V^2])\\&=&
(\tilde{\nabla}_{V^1}\phi)(V^2)+\phi(\tilde{\nabla}_{V^1})V^2-
(\tilde{\nabla}_{V^2}\phi)(V^1)-\phi(\tilde{\nabla}_{V^2}V^1)-\phi([V^1,V^2])\\
&=&\tilde{\nabla}_{V^1}\phi(V^2)-\tilde{\nabla}_{V^2}\phi(V^1)\\
&=&{1\over 2} d\phi(V^1, V^2)
\end{eqnarray*}
by the standard formula, as the pointwise connection $\tilde \nabla$ has no torsion.
\hfill\rule{2mm}{2mm}

\end{document}